\documentclass[12pt]{amsart}
\newcommand{\kcalomega}{\mathcal{K}_{[\omega_0]}}
\def\XXint#1#2#3{{\setbox0=\hbox{$#1{#2#3}{\int}$ }
\vcenter{\hbox{$#2#3$ }}\kern-.6\wd0}}

\renewcommand{\d}{\partial}
\newcommand{\sym}{{\operatorname{Sym}}}

\usepackage{graphicx}

\usepackage{amsmath, amssymb,latexsym}
\usepackage{verbatim}
\usepackage{mathtools}
\usepackage{color}
\usepackage[active]{srcltx}
\usepackage{graphicx}
\usepackage{hyperref}

\newcommand{\K}{{\mathbf K}}

\def\p{\partial}
\def\bp{\bar\partial}

\def\d{\delta}

\usepackage{color}

\newcommand{\crit}{{\operatorname {Crit}}}
\newcommand{\beq}{\begin{equation}}
\newcommand{\eeq}{\end{equation}}

\newcommand{\dist}{{\operatorname{dist}}}
\newcommand{\Li}{{\operatorname{Li}}}

\usepackage[active]{srcltx}

\newenvironment{example}{\medskip\noindent{\it Example:\/} }{\medskip}

\newcommand{\Var}{{\bf{Var}}}
\newcommand{\var}{{\operatorname{Var}}}

\newcommand{\I}{{\mathbf I}}

\newcommand{\szego}{Szeg\H o }

\newcommand{\inv}{^{-1}}
\newcommand{\kahler}{K\"ahler }
\newcommand{\sqrtk}{\sqrt{k}}
\newcommand{\wt}{\widetilde}

\newcommand{\PP}{{\mathbb P}}

\newcommand{\R}{{\mathbb R}}
\newcommand{\C}{{\mathbb C}}

\newcommand{\Z}{{\mathbb Z}}

\newcommand{\CP}{\C\PP}

\renewcommand{\d}{\partial}
\newcommand{\dbar}{\bar\partial}
\newcommand{\ddbar}{\partial\dbar}

\newcommand{\E}{{\mathbf E}}

\renewcommand{\phi}{\varphi}

\newcommand{\bcal}{\mathcal{B}}
\newcommand{\ccal}{\mathcal{C}}
\newcommand{\dcal}{\mathcal{D}}
\newcommand{\ecal}{\mathcal{E}}

\newcommand{\hcal}{\mathcal{H}}

\newcommand{\lcal}{\mathcal{L}}

\newcommand{\ncal}{\mathcal{N}}
\newcommand{\ocal}{\mathcal{O}}
\newcommand{\pcal}{\mathcal{P}}

\newcommand{\scal}{\mathcal{S}}

\newcommand{\al}{\alpha}
\newcommand{\be}{\beta}
\newcommand{\ga}{\gamma}

\newcommand{\la}{\lambda}
\newcommand{\ep}{\varepsilon}
\newcommand{\de}{\delta}
\newcommand{\De}{\Delta}
\newcommand{\om}{\omega}
\newcommand{\Om}{\Omega}

\newcommand{\half}{{\frac{1}{2}}}
\newcommand{\vol}{{\operatorname{Vol}}}

\newcommand{\SU}{{\operatorname{SU}}}
\newcommand{\SL}{{\operatorname{SL}}}

\newcommand{\FS}{{{\operatorname{FS}}}}
\renewcommand{\phi}{\varphi}

\newtheorem{theo}{{\sc Theorem}}[section]
\newtheorem{cor}[theo]{{\sc Corollary}}

\newtheorem{defin}[theo]{{\sc Definition}}

\newtheorem{conj}[theo]{{\sc Conjecture}}
\newtheorem{lem}[theo]{{\sc Lemma}}

\newtheorem{prop}[theo]{{\sc Proposition}}

\title[Stochastic K\"ahler Geometry]{Stochastic K\"ahler Geometry: From random zeros to  random metrics}

\author{Bernard Shiffman }
\address{Department of Mathematics,  Johns Hopkins University, Baltimore, MD 21218, USA }
\email{bshiffman@jhu.edu}

\author{Steve Zelditch}  
\address{Department of Mathematics, Northwestern  University, Evanston, IL 60208, USA}

\thanks{Research of the second author partially supported by    NSF grant DMS-1810747}

\begin{document}

\begin{abstract} We provide a survey of results on the statistics of random sections of holomorphic line bundles on \kahler manifolds, with an emphasis on the resulting asymptotics  when a line bundle is raised to increasing tensor powers. We conclude with a brief discussion of the `Bergman' \kahler metrics induced by these random sections.
\end{abstract}

\maketitle

\tableofcontents

\section*{Introduction}

  Stochastic \kahler geometry  refers to the study of probabilistic problems in complex algebraic or analytic geometry in the setting of K\"ahler manifolds $(M^m, J, \omega)$ of any complex dimension $m $.   It concerns  random fields  
on a \kahler manifold  which are defined in terms of  the complex structure $J$
and \kahler form $\omega$. The basic random fields are holomorphic
sections $s \in H^0(M, L^k)$ of powers of a holomorphic Hermitian line bundle
$(L, h) \to (M, \omega)$. From these holomorphic fields one can construct 
random  complex submanifolds $Z_{\vec s}$  (zero sets of one or several sections), random embeddings into complex projective spaces $\CP^N$, and random `Bergman' or `Fubini-Study' \kahler metrics induced by the embeddings.  Zero sets and embeddings both determine positive (1,1)-forms
$\omega = {i} \ddbar \phi$, where $\phi$ is a  psh (plurisubharmonic)  function.  Although the behavior of random zero sets in the high tensor power limit $k \to \infty$ is the heart of stochastic K\"ahler geometry, the same techniques often apply with little change to random K\"ahler metrics and other more general objects.  The goal of this survey is to review some of the main results on random zero sets and also to briefly discuss these generalizations to random K\"ahler metrics. 

Most results of stochastic \kahler geometry to date pertain to the asymptotics of probabilistic invariants such as  distribution and correlation functions of zeros and of critical points  
as the degree $k \to \infty$.  One of the main results is universality of the limit of rescaled invariants on small balls of radius $k^{-\half}$. Recently such scaling limits have been used to study the local topology of random zero sets. Another focal point is on  the asymptotic normality of linear statistics, showing that fluctuations of linear statistics, i.e.\  integrals  $\int_{Z_s} \psi $  of a test form $\psi$  over  the zeros of random sections,   tend to Gaussian random variables determined by the variance current. Asymptotic normality of integrals against a  random positive $(1,1)$ form $\omega$ is equivalent to 
asymptotic normality of the potential $u$ of $\omega$, and that is the way it is often stated in the physics literature (e.g. \cite{CLW,CLW2}). We call attention to  some natural ensembles of potentials for which asymptotic normality is as yet unknown: linear statistics for critical points and for zeros of codimension greater than 1.

As the reference to potentials indicates, the  unifying theme is that of random psh functions. 
A \kahler metric is defined as a mixed Hessian
 $\omega = i \ddbar u$ of  a local psh  function  $u$, known as the  `\kahler potential'. Zero sets are also defined as $Z_{\vec f}: = (\frac i{2\pi}\,\ddbar u)^q$
 where $u = \log \sum_{j = 1}^q |f_j(z)|^2$ with $q \leq m = \dim_{\C} M$ and $f_j$ are local holomorphic functions. The same formula when $q > m$ is a way to define a 
 smooth \kahler metric, and $Z_{\vec f}$ can be viewed as a `singular \kahler metric'. As this suggests, many results about random zero sets have analogues for random smooth metrics. If $q = d_k: = \dim H^0(M, L^k)$
 and if $\vec f$ is a basis of $H^0(M, L^k)$ then  $i \ddbar \log \sum_{j = 1}^{d_k} |f_j(z)|^2$  is known as a Bergman metric of degree $k$.

 One  of the themes of stochastic \kahler geometry is the response of the probabilistic results to changes in the input geometry.
 By `geometry' we mean line bundles $L \to M$, Hermitian
metrics $h$ on $L$, curvature forms $\Theta_h$ and the `quantization' of
Hermitian metrics $h$ (together with a choice of measure $\nu$ on $M$) as inner products $G(h^k, \nu)$ on spaces
$H^0(M, L^k)$ of holomorphic sections of powers of $L$. The inner product
determines a Gaussian measure $\gamma_{h^k, \nu}$ on $H^0(M, L^k)$, and this  provides  the  notion of `random polynomial' or more generally
`random section'.  The geometric
language is  useful (and even necessary) to formulate generalizations
of logarithmic potential theory and random polynomial theory on $\C$ to
compact Riemann surfaces or to higher dimensional complex manifolds.
Holomorphic sections of line bundles are the analogues on a compact manifold $M$ of holomorphic
functions  on $\C^n$, and specifically $H^0(M, L^k)$ is the replacement
for polynomials of degree $\le k$.  

We may contrast stochastic \kahler geometry with the much-studied
one-dimensional theory of stochastic (or, random) conformal geometry. 
Conformal stochastic geometry is a highly developed field of probability, mathematical physics and complex analysis. It contains such subfields as  SLE, the  quantum Hall effect, Hele-Shale flow, and Liouville quantum gravity in mathematical physics, and probabilistic problems in one  dimensional complex analysis.  As the name `conformal' suggests, it is strictly a complex one-dimensional
theory.  
 The key difference is that stochastic conformal geometry is concerned with conformally-invariant ensembles of real objects such as the  Gaussian free field (GFF),  SLE curves,  Coulomb gas point processes,  or random LQG area forms in Liouville quantum gravity    \cite{AHM11, Du06,KN13}. The key objects are often random fractals. In stochastic \kahler geometry, the emphasis is on holomorphic fields and the  
 objects they induce in complex geometry.
  
In this survey, we only refer briefly to results in the complex one dimensional case, although it is a very rich field. Moreover, many of the recent constructions on higher dimensiional K\"ahler manifolds use ideas that originated in the probabilistic study of real algebraic manifolds and zero sets of random real functions, in particular ideas 
stemming from the work of Nazarov--Sodin \cite{NS09} on counting connected components of spherical harmonics and other random real functions,  and their  topological 
applications due to Sarnak--Wigman \cite{SW19},  Canzani--Sarnak \cite{CS19} and others on Betti numbers and combinatorial configurations.  We omit these important results because they would take us too far afield.
\section{\label{BACKGROUND} Background}

In this section we introduce some background and notation pertaining to random holomorphic sections of positive Hermitian line bundles. 

Let $(M,L)$ be an $m$-dimensional compact complex manifold polarized  with a Hermitian holomorphic  line bundle $(L,h)$. We
consider a local holomorphic frame $e_L$ over a trivializing chart
$U$. If  $s = f e_L$ is a holomorphic section of $L$ over $U$,  its Hermitian norm is
given by
$\|s(z)\|_h = e^{-\phi_h}|f(z)|$ where \begin{equation}
\label{a} \phi_h(z): = -\log\|e_L(z)\|_h\;. \end{equation} 

The curvature form of
$(L,h)$ is given locally by
$\Theta_h= 2\ddbar \phi_h$, and the
{\it Chern form\/} $c_1(L,h)$ is given by
\begin{equation}\label{chern}c_1(L,h)=\frac{\sqrt{-1}}{2 \pi}
\Theta_h=\frac{\sqrt{-1}}{ \pi}\,\d\dbar\phi_h\;.\end{equation} 
We now assume that the Hermitian metric $h$ has strictly positive
curvature and we give
$M$ the \kahler form \begin{equation}\label{omega}\om_h:= i\ddbar\phi_h=\pi
c_1(L,h)\;.\end{equation}

\subsection{\label{inner+gaussian}From metrics and measures to inner products and Gaussian measures}

We denote by $H^0(M, L^{k})$ the space of global holomorphic
sections of
$L^k=L\otimes\cdots\otimes L$.   The metric $h$ induces Hermitian
metrics
$h^k$ on $L^k$ given by $\|s^{\otimes k}\|_{h^k}=\|s\|_h^k$. 

We let $\nu$ denote a (finite, positive) Borel measure on $M$. Together, the data $(h, \nu)$  induces Hermitian inner products $G(h^k, \nu)$ on the spaces  $H^0(M,
L^k)$ of global holomorphic sections of powers $L^k \to M$ given by
\begin{equation}\label{inner}\langle s_1, \overline{s_2} \rangle_k = \langle s_1, \overline{s_2} \rangle_{G(h^k, \nu)} : = \int_M\langle s_1(z) ,\overline{s_2(z)}\rangle_{h^k} \,  d\nu(z), \quad s_1,s_2\in H^0(M,L^k)\,.\eeq

In turn, each inner product on $H^0(M, L^k)$ induces an 
orthonormal basis $\{S_1^k,\dots,S_{d_k}^k\}$ and  associated Gaussian
measure $\gamma_{h^k,\nu}$ given  by the formula,
\begin{equation}\label{gaussian}d\gamma_{h^k,\nu}(s^k):=\frac{1}{\pi^{d_k}}e^
{-|c|^2}dc\,,\quad s^k=\sum_{j=1}^{d_k}c_jS^k_j\,,\quad
c=(c_1,\dots,c_{d_k})\in\C^{d_k}\,,\end{equation} where
 $dc$ denotes $2d_k$-dimensional
Lebesgue measure.   The measure $\gamma_{h^k,\nu}$ is characterized by the
property that the $2d_k$ real variables $\Re c_j, \Im c_j$
($j=1,\dots,d_k$) are independent Gaussian random variables with
mean 0 and variance $1/2$; equivalently,
$$\E(c_j) = 0 = \E(c_j c_l),\;\;\; \E(c_j \bar{c}_l) =
\delta_{jl}.$$ Here, $\E$ denotes the expectation.

The inner product $G(h^k, \nu)$ further induces  an associated spherical
measure  on the unit sphere $SH^0(M, L^{k})$ in
$H^0(M, L^{k})$ with respect to $G(h^k, \nu)$.
In this survey, we restrict our discussion to inner products where $\nu$ is the volume form of $M$. For  results with more general metrics and measures, see for example \cite{Ber09, BS07, BBW}.

\subsection{Polynomials and holomorphic sections of line bundles}

The space  ${\it Poly}_k$ of univariate polynomials of degree $k$ is a complex vector space of dimension $k + 1$. 
The `$\SU(2)$  inner product' on ${\it Poly}_k$  may be written in the form
$$\langle f_1,\bar f_2\rangle=\frac i2\int_{\C} f_1(z) \overline{f_2(z)} e^{- k \log(1 + |z|^2)} \frac{dz \wedge d\bar{z}}{(1 + |z|^2)^2}=\frac i2\int_{\C} f(_1z) \overline{f_2(z)}  \frac{dz \wedge d\bar{z}}{(1 + |z|^2)^{k+2}}\,,$$ $f_1,f_2\in{\it Poly}_k$.

This has a simple geometric interpretation: namely, we view polynomials of degree $k$ as
holomorphic sections of the line bundle $\ocal(k)=L^k$,  where $L\to \CP^1$ is the hyperplane section bundle. We give $L$ the Hermitian metric
$\|e_
L\|_h= e^{-\frac 12 \log(1 + |z|^2)}$, where $e_L=1$. Then 
$\om_h=\frac i2 \frac{dz \wedge d\bar{z}}{(1 + |z|^2)^2}$ is the usual area form
on $\hat{\C} = \CP^1$, and  $\langle f_1,\bar f_2\rangle=\langle f_1,\bar f_2\rangle_{G(h^k,\om_h)}$.

 For multivariable polynomials, we let $ M=\CP^m$ with $L\to\CP^m$ the hyperplane section bundle $\ocal(1)$, so that the space of global sections $H^0(\CP^m,L)$ consists of the linear functions $f(z)=\sum_{j=0}^mc_jz_j$ on $\C^{m+1}$. Then $H^0(\CP^m,L^k)$ is the vector space ${\it Poly}_k^m$ of homogeneous polynomials of degree $k$ on $\C^{m+1}$, which we identify with the space of polynomials of degree $\leq k$ in the variables $z_1,\dots,z_m$ by setting $z_0=1$. If we let $w_j=z_j/z_0$,\;  $1\le j\le m$, be local coordinates on $\CP^m$ and we give $L=\ocal(1)$ the  Hermitian metric  $\|e_L\|_h=(1+\|w\|^2)^{-1/2}$, then $\phi_h=\frac 12\log(1+\|w\|^2)$ and $\om_h= \frac i2\ddbar\log(1+\|w\|^2)$, the Fubini-Study metric on $\CP^m$. Then the volume form $$dV=\frac 1{m!}\om_h^m = (1+\|w\|^2)^{-m-1}\,d_{2m}w\,,$$ where $d_{2m}w$ is Euclidean volume. We then have the $\SU(m+1)$-invariant inner product
 $$\langle f_1,\bar f_2\rangle_{G(h^k,dV)}=\int_{\C^m} \frac{ f_1(z)\, \overline{f_2(z)} }{(1 + \|z\|^2)^{k+m+1}}\,d_{2m}z\,,\quad f_1,f_2\in{\it Poly}_k^m\;.$$

\subsection{Asymptotics of Bergman kernels on positive line 
bundles}

We let $\lcal^{2}(M,L^{k})$ denote the  $\lcal^2$  sections of $L^k\to M$ with respect to the inner product $G(h^k,dV)$, where $dV=\frac 1{m!}\om_h^m$.

We define the  {\it Bergman kernel} as the orthogonal projection $$B_k(z,w): \lcal^{2}(M,L^{k})\rightarrow H^{0}(M,L^{k})\,.$$ Then 
\beq \label{bergmankernel} B_k(z,w) = \sum_{j=1}^{d_k} S_j^k(z) \otimes \overline{S_j^k(w)}\,,\eeq
{where $\{S_1^k,\cdots, S_{d_k}^k\}$ is an orthonormal basis of $H^0(M,L^k)$ with respect to  $G({h^k},dV)$.}
Along the diagonal, the contraction of the Bergman kernel  is 
\beq\|B_k (z,z)\| =\sum_{j=1}^{d_k} \|S^k_j(z)\|_{h^k}^2\,.\eeq
In the case where the curvature form $\Theta_h$ of the Hermitian line bundle $(L,h)$ is everywhere positive, we  have the following {\it Tian--Yau--Zelditch  asymptotic expansion} \cite{Ca97, Z97}:
\begin{equation}\label{seg} \|B_k (z,z)\|\sim \frac 1{\pi^m}[k^m+a_1(z)k^{m-1}+a_2(z)k^{m-2}+     \cdots]\,,
\end{equation} where each coefficient
$a_j (z)$ is a polynomial of the curvature and its covariant derivatives. Formulas for the first three coefficients were given by Lu \cite{Lu00}. In particular, $a_1(z)$ equals one-half the scalar curvature of $\omega_h$. 

\begin{example}In the case of $(\CP^m, \omega_{FS})$ with the line bundle $(\mathcal O(1), h_{FS})$, the Bergman kernel is easily computed to be a constant along the diagonal \cite{BSZ00}:
\begin{equation} \label{cp}\|B_k (z,z)\|=\frac{(k+m)!}{\pi^m\,k!}\,.
\end{equation}
\end{example}

Recall that  the inner product $G(h^k,dV)$ induces the Gaussian field $(H^0(M,L^k),\gamma_{h^k,dV})$, where $\gamma_{h^k,dV}$ is given by \ref{gaussian}. { In fact, the Bergman kernel $B_k(z,w)$ can be interpreted as the covariance function for the Gaussian field $\big(H^0(M,L^k),\,\gamma_{h^k,dV}\big)$:}
\beq\label{2point} \E\big(s^k(z)\otimes\overline{s^k(w)}\big) = B_k(z,w)\,,\eeq
where $\E$ denotes the expected value with respect to $\gamma_{h^k,dV}$.

\begin{proof} Apply $ \E( c_j \bar c_l) =
\de_{jl}$ to
$$ \E \left(s^k(z)\otimes \overline{s^k(w)}\right) = \E\left(\sum_{j=1}^{d_k}c_jS^k_j(z)\otimes\overline{\sum_{l=1}^{n}c_lS^k_l(w)}\right). $$ \end{proof}

\subsection{Off-diagonal scaling asymptotics  of the \szego kernel}\label{offdiag} To provide asymptotics for the Bergman kernel off the diagonal, it is convenient to lift the Bergman kernel to the circle bundle $X$ of the dual bundle to $L$. To describe the lifted kernel, we let $L^*\to M$ denote the dual line bundle to $L\to M$ with the dual metric $h^*$, and we let $X:=\{\lambda\in L^*:\|\la\|_{h^*}=1\}$.  We regard a section $s^k\in H^0(M,L^k)$ as a function on $X$ by setting $$s^k(\lambda) = (\la\otimes \cdots\otimes \la,s^k(z)),\quad \lambda\in L^*_z\,.$$ and we note that $s^k$ is $k$-equivariant: $s^k(e^{i\theta}\la)=e^{ik\theta}s^k(\la)$. We assume that $(L,h)$ has  positive curvature; then $X$ is the boundary of the strictly pseudoconvex disk bundle $\{\la\in L^*: \ell(\la)<1\}$ where  $\ell(\la)=\|\la\|^2_{h^*}$. We let $\Pi:\lcal^2(X)\to \hcal^2(X)$ denote the orthogonal projection to the space { $\hcal^2(X)$ of square-integrable} CR functions on $X$, where we give $X$ the volume form $$\frac i{2\pi\, m!}\dbar \ell\wedge (i\ddbar\ell)^m= \frac i{2\pi}\dbar \ell\wedge dV_M\,.$$ Then $\Pi=\bigoplus_{k=0}^\infty\Pi_k$, where $\Pi_k:\lcal^2(X)\to \hcal^2_k(X)$ is the orthogonal projection onto  the space of $k$-equivariant functions $\hcal^2_k(X)$ in $\hcal^2(X)$. Indeed, $\hcal^2(X)=\bigoplus_{k=0}^\infty\hcal^2_k(X)$, where   $\hcal^2_k(X)\approx H^0(M,L^k)$.
We call $\Pi_k(x,y)$ the ($k$-th) \szego kernel; the sum $\Pi(x,y)$  is the classical \szego kernel for the strictly pseudoconvex boundary $X$. 

To relate the Bergman kernel $B_k(z,w)$ to the \szego  kernel $\Pi_k(z,\theta;w,\phi)$, we use  a local frame $e_L$ to write $S^k_j=F_j^ke_L^{\otimes k}$.
Recalling \eqref{bergmankernel}, we have 
\begin{eqnarray*}B_k(z,w)&=&\textstyle \left( \sum_{j=1}^{d_k}F_j^k(z)\overline{F_j^k(w)}\right )e_L(z)^{\otimes k}\otimes \overline{e_L(w)}^{\otimes k}\;, \\ \Pi_k(z,\theta_1;w,\theta_2)&=&\textstyle e^{ik(\theta_1-\theta_2)}e^{-k\phi(z)-k\phi(w)}\sum_{j=1}^{d_k}F_j^k(z)\overline{F_j^k(w)}\,.\end{eqnarray*}
Here, $(z,\theta)$ denotes the point
$e^{i\theta}\|e_L(z)\|_he^*_L(z)\in X$. Thus
\begin{eqnarray*} |\Pi_k(z,\theta_1;w,\theta_2)|= \|B_k(z,w) \|
\;,\\ \Pi_k(z,z):=\Pi_k(z,0;z,0)= \|B_k(z,z) \|\;.\end{eqnarray*} 

The asymptotics of the Bergman kernal are used in Section \ref{section-zeros} to study the distributions of zeros of a random section $s^k\in H^0(M,L^k)$. In particular, the off-diagonal asymptotics of the Bergman kernel provides information on correlations and variances of random zeros. To this end, a  general asymptotic expansion was given in \cite{ShZ02} and further clarified in \cite{ShZ08} as follows:

\begin{theo}
\label{near} Let $(L,h)\to
(M,\om_h)$ be a positive holomorphic line bundle over a compact \kahler manifold.  Let $z_0\in M$ and choose local coordinates
$\{z^j\}$ in a neighborhood of $z_0$ so that $z_0=0$ and
$\Theta_h(z_0)=\sum dz^j\wedge d\bar z^j$.  Then
\begin{multline*}\frac{\pi^m}{k^m}\,\Pi_k\left(\frac{z}{\sqrtk},\frac{\theta_1}{k };
\frac{w}{\sqrtk},\frac{\theta_2}{k}\right)=e^{i (\theta_1-\theta_2)+i\Im
(z\cdot \bar w)-\half |z-w|^2}\\ \times \left[1+ \sum_{r = 1}^{n} k^{-r/2}
p_{r}(z,w) + k^{-(n +1)/2} R_{kn}(z,w)\right]\;,\end{multline*}
where $p_r$ is a polynomial in $(z,\bar z,w,\bar w)$ of the same parity as
$r$, and
$$|\nabla^jR_{kn}(z,w)|\le C_{jn\ep b}k^{\ep}\quad \mbox{for }\
|z|+|w|<b\sqrt{\log k}\,,$$ for
$\ep,b\in\R^+$,  $j,k\ge 0$. Furthermore, the constant  $C_{jk\ep
b}$ can be chosen independently of $z_0$. 
\end{theo}
Here, $\nabla^j$ stands for the $j$-th covariant derivative.

The theorem shows that on K\"ahler manifolds, there is a characteristic length scale associated to the $k$-th power 
$L^k \to M$ of a positive line bundle: the Planck scale $\frac{1}{\sqrt{k}}$.  It arises in the following ways:

\begin{itemize}

\item The \szego kernel $\Pi_{k}(z,z_0)$ is of size $ \simeq k^m$  for $ \dist(z,z_0)< \frac{b}{\sqrt{k}}$, and then decays rapidly outside the ball. 

\item On the length scale $\frac{1}{\sqrt{k}}$, all \kahler manifolds and positive  line bundles look alike {in the scaling limit:} they
all look like the (trivial) line bundle { $\C^{m + 1} \to \C^m$ with the Euclidean \kahler form on $\C^m$.} 

\item Correlations become universal on this length scale.

\end{itemize}

Specific formulas for the coefficients in the off-diagonal expansion using Bochner coordinates are given by Lu--Shiffman in \cite{LSh}. For real-analytic metrics, \cite{HLX20} gives symptotics on an enlarged length scale. 
A detailed study of the off-diagonal asymptotics is given in the book of Ma--Marinescu \cite{MaMa07} using different techniques  involving normal coordinates instead of holomorphic coordinates.

Away from the diagonal, we have the following decay estimate \cite{ShZ08}:
\begin{theo}\label{far} Let $(L,h)\to
(M,\om_h)$ be as above. For
$b>\sqrt{j+2\al+2m}\,$, $j,\al\ge 0$, we have
$$ \nabla^j
\Pi_k(z,0;w,0)=O(k^{-\al})\qquad \mbox{uniformly for }\ \dist(z,w)\ge
b\,\sqrt{\frac {\log k}{k}} \;.$$
\end{theo}

 In particular, our variance formulas are
expressed in terms of the {\it normalized Bergman kernel}
\begin{equation}\label{Pk} P_k(z,w):=
\frac{\|B_k(z,w)\|}{\|B_k(z,z)\|^\frac 12\,\| B_k(w,w)\|^\frac
12}\,.\end{equation}
which is the square root of the so-called {\it Berezin kernel}.  Note that $0\le P_k(z,w)\le 1$ by Cauchy-Schwarz, and $P_k(z,z)=1$.

Theorems \ref{near}--\ref{far} yield the following counterparts for the normalized kernel $P_k(z,w)$:

\begin{prop} \label{DPdecay} Let $(L,h)\to
(M,\om_h)$ be a positive holomorphic line bundle over a compact \kahler manifold. For
$b>\sqrt{j+2\al}$,
$j,\al\ge 0$, the normalized Bergman kernel satisfies the asymptotic estimate
$$ \nabla^j
P_k(z,w)=O(k^{-\al})\qquad \mbox{uniformly for }\ d(z,w)\ge
b\,\sqrt{\frac {\log k}{k}} \;.$$
\end{prop}

\begin{prop} \label{better} Using the hypotheses and notation of Theorem \ref{near}, we have the following asymptotics for the normalized Bergman kernel near the diagonal: 

For $\ep,b>0$, there are constants
$C_j=C_j({M,\ep,b})$, $j \ge 2$, independent of the point $z_0$, such that
$$  P_k\left(\frac z{\sqrtk},\frac w{\sqrtk}\right) =
e^{-\frac 12 |z-w|^2}[1 + R_k(z,w)]\;,$$ where
$$ \begin{array}{c}|R_k(z,w)|\le \frac {C_2}2\,|z-w|^2k^{-1/2+\ep}\,, \quad
|\nabla R_k(z)|
\le C_2\,|z-w|\,k^{-1/2+\ep}\,,
\\[8pt] |\nabla^jR_k(z,w)|\le C_j\,k^{-1/2+\ep}\quad j\ge 2\,,\end{array}$$
for $|z|+|w|<b\sqrt{\log k}$.\end{prop}

\section{Random zero sets}\label{section-zeros}

 {  

We now consider zero sets
 $$Z_s = \{z \in M: s(z) = 0\} $$ of Gaussian random  holomorphic sections $s\in H^0(M,L)$. In  the case where $M$ is a compact Riemann surface $C$ (complex dimension 1), 
 the zero set  $Z_s$  is a finite set $\{\zeta_1, \dots, \zeta_k
 \}$ of points in $C$.  For example, if $C=\CP^1$ and $L=\ocal(k)$, then $s$ is a polynomial on $\C\subset\CP^1$ of degree $\le k$ and the zero set consists of the roots of $s$ (and the point at infinity if $\deg s <k$).
 
From the probabilistic viewpoint, the zeros of a random
holomorphic section define a {\it point process} on
$C$, that is, a measure on the configuration space ${\rm Conf}(C)$ of
finite subsets of $C$ (where the points may have positive integral multiplicities). Each holomorphic section gives rise to the
discrete set of its zeros, and the  point
process is the probability measure on ${\rm Conf}(C)$ induced by the  probability measure on the vector space $H^0(C,L)$. A probability measure on ${\rm Conf}(C)$ is determined by
its $n$-point correlations $\K_{n} (z_1,
\dots, z_n)$, $n\ge 1$, which are the probability densities (in $C^n=C\times\cdots\times C$) that $z_1, \dots, z_n \in C$ are the (simultaneous) zeros of a random section.  For example, the {\it pair correlation} $\K_2(z_1,z_2)$ determines whether the zeros tend to cluster or to `repel' each other.

The zero set $s^{-1}(0)=\{\zeta_1,\dots,\zeta_k\}\in {\rm Conf}(C)$ of zeros of a section $s$ yields the {\it normalized  empirical measure} 
$$\frac{1}{k} Z_s = \frac{1}{k} \sum_j \delta_{\zeta_j}, $$
(again, counting multiplicities), so that the point  process can be considered as a measure on the space of probability measures on $C$ with discrete support.
}
 Here,  $\delta_z$ is the Dirac delta-function at $z$. Thus the { normalized}
 empirical measure of zeros,
$$(\frac{1}{k} Z_s,\psi)= \frac{1}{k} \sum_j \psi(\zeta_j), { \quad \psi\in\ccal(C)},$$
is a random probability measure on $C$. Its expectation is a
measure called the {\it expected distribution of zeros.}  

For $m=\dim M \ge 2$, $Z_s$ is the {\it current of integration} over the zero set of $s$:
\beq\label{current}  (Z_s,\psi) = \int_{Z_s'}\psi,\qquad \psi \in\dcal^{m-1,m-1}(M)\,,\eeq
where $Z_s'$ is the set of smooth points (counted with multiplicities) of the analytic hypersurface $\{\zeta:f(\zeta)=0\}$. In Section \ref{ppm}, we discuss point processes of simultaneous zeros of $m$ holomorphic sections on $M$.

In \cite{ShZ99}, we showed the following:

\begin{theo} \label{EZsk} Let $(L, h) \to  (M,\om_h)$ be a positive line bundle over a compact \kahler manifold. Then   $$ \frac 1k\E (Z_{s^k})\to  \frac1\pi \,\omega_h  $$   weakly in the sense
of measures, where $\E$ is the expectation with respect to the Gaussian measure  $ \ga_{h^k,dV}$  on $H^0(M,L^k)$.  In fact,
\begin{equation}\label{mean}\frac 1 {k}\,\E ( Z_{s^k}, \psi ) =\frac 1\pi
\int_M  \omega_h \wedge\psi\ +\ O\left(\frac 1{k^2}\right)\,,\qquad
\psi\in\dcal^{m-1,m-1}(M)\;.\end{equation}   Thus, if $\{s^k\in H^0(M,L^k)\}$ is a sequence of independent random sections, then
$$\frac 1kZ_{s^k} \to \frac 1\pi\,\omega_h\quad a.s.$$
\end{theo}

Precisely, we form the probability space $\scal:=\prod_{k = 1}^{\infty} H^0(M, L^k)$
with the product measure. Its elements are sequences $\{s^k\}$
of independent random sections. (In \cite{ShZ99}, \eqref{mean} was stated with
remainder term $O(\frac 1k)$ in place of  $O(\frac 1{k^2})$\,.)

In
particular, we have
$$\lim_{k \rightarrow \infty}\frac{1}{k} \E\vol_{2m-2}(Z_{s^k}\cap U)\to \frac m\pi \vol_{2m}(U)\,,$$ for $U$ open in $M$.
In the Riemann surface case ($m=1$), 
$$\lim_{k \rightarrow \infty}\frac{1}{k} \E \#\{Z_{s^k}\cap U\}={\frac 1\pi\operatorname{Area}(U)}\;.$$
We outline the proof of Theorem \ref{EZsk} in the next section.

\subsection{Poincar\'e-Lelong formula}

 In complex dimension one, if $f(z)$ is a holomorphic
function on a domain in $\C$, then the fundamental solution $\Delta \log |z|^2 = 2 \pi \delta_0$ of the Laplace operator immediately yields
\begin{equation}
\label{PL1}Z_f  = \sum_{f(\zeta)=0} \delta_{\zeta} =  \frac{i}{2 \pi} \ddbar \log |f|^2  = \frac{i}{2\pi}
\frac{\partial^2 \log |f|^2}{\partial z \partial \bar{z}} dz \wedge d \bar{z}\,, \end{equation}
as a singular $(1,1)$-current. 

In higher dimensions, we similarly have {(see\cite{Lelong})}\beq
\label{PLf}  Z_f  = \frac{i}{2 \pi} \ddbar \log |f|^2 \;,\end{equation}
{where $Z_f \in\dcal^{\prime 1,1}(M)$ denotes the current of integration given in \eqref{current}.}

 For a section ${s^k}= f e_L^{\otimes k}\in H^0(M, L^k)$ of a Hermitian holomorphic line bundle $L\to M$, we then have by\eqref{PLf} the {\it Poincare-Lelong formula,}
\begin{equation}
\label{PL}  Z_{s^k}  = \frac{i}{ \pi} \ddbar \log |f| = \frac{i}{\pi}
\ddbar \log \|{s^k}\|_{h^k} + \frac k\pi\omega_h\;.\end{equation}

{Averaging \eqref{PL}, we obtain:

\begin{theo}\label{EPL} Let $\{S_j^k\}$ be an orthonormal basis of $ H^0(M,L^k)$. 
Write $S_j^k =f_j^k e_L^{\otimes k}$.  Then,
$$\frac 1k\E ( Z_{s^k}) = \frac{\sqrt{-1}}{2 \pi k } \ddbar \log \sum_{j=1}^{d_k}
|f^k_j|^2= \frac{\sqrt{-1}}{2 \pi k } \ddbar \log
\|B_{k}(z,z)\| +\frac 1\pi \om,
$$
where{ we recall that $B_{k}$} is the Bergman kernel. 

\end{theo}
\begin{proof}

Let $s = \sum_j a_j S_j^k$ and write it as $s=\langle \vec a, \vec S^k \rangle
= \langle \vec a, \vec f \rangle e^k_L$. 
 Let $\psi \in \dcal^{m-1,m-1}(M)$. Then
 $$ \E \langle \frac{1}{k} [Z_s^k]), \psi \rangle = \frac{\sqrt{-1}}{\pi k} \int_{\C^{d_k}} d\gamma_k(a)\int_M \ddbar \log
|\langle \vec a, \vec f\rangle|\wedge \psi\;.$$ To compute the integral, we write $\vec f = |\vec f|
\vec u$ where $|\vec u| \equiv 1.$ Evidently, $\log |\langle \vec a, \vec f\rangle| =
\log |\vec f| + \log |\langle \vec a, \vec u \rangle|$. The first term gives
\begin{equation} \frac{\sqrt{-1}}{\pi k} \int_M \ddbar \log
|\vec f| \wedge  \psi = \frac{\sqrt{-1}}{2\pi k} \int_M \ddbar \log \|B_k(z,z)\|\wedge
  \psi  + \frac 1\pi \int_M \omega\wedge \psi. \end{equation}

 We
now look at the second term.  We have \begin{multline}\label{term=0}
\frac{\sqrt{-1}}{\pi}\int_{\C^{d_k}} d\gamma_k(a) \int_M \ddbar \log
|\langle \vec a, \vec u\rangle| \wedge \psi \\=
\frac{\sqrt{-1}}{ \pi}\int_M  \ddbar \left[\int_{\C^{d_k}} \log |\langle \vec a, \vec u \rangle| \,d\gamma_k(a)
 \right] \wedge \psi =0,\end{multline}
since the average $\int \log |\langle \vec a, \vec u \rangle|
d\gamma_k(a)$ is a constant independent of $\vec u$ for $|\vec u|=1$, and thus
the operator $ \ddbar$ kills it. \end{proof}

Combining Lemma \ref{EPL} with the Bergman kernel asymptotics \eqref{seg} yields Theorem \ref{EZsk}.

}

\subsection{\label{PCFZSECT}Correlation  of zeros}

In this section, we  discuss  $n$-point `correlations' between
zeros, or `joint intensities', of   random sections $s^k\in H^0(M,L^k)$ of powers of a positive line bundle. We first consider pair correlations ($n=2$):
the pair correlation current for random zeros is defined by 
 \begin{equation} \label{Kk} \K_2^k(z,w):=\E\left(Z_{s^k}(z)\otimes Z_{s^k}(w)\right);\end{equation}
i.e.,  for  test forms $ \psi_1, \  \psi_2 \in \dcal^{m-1,m-1}(M)$,
\begin{equation}\left(\K_2^k(z,w), \psi_1(z)\otimes \psi_2(w)\right) :=\E\left[\left(Z_{s^k}, \psi_1\right)\otimes\left(Z_{s^k}, \psi_2\right)\right].\end{equation}

In the case of complex dimension 1, the zeros form a point process, as discussed above, and the pair  correlation measures take the form
$$   \K_2^k(z,w) = [\Delta]\wedge (  \K_{1}^k(z)
\otimes 1)  + \tilde K_2^k(z,w)\,\K_1^k(z)\otimes \K_1^k(w) \,,$$ where $[\Delta]$ denotes the current of integration along the
diagonal $\Delta=\{(z,z)\}\subset C\times C$, and $\K^k_1=\E(Z_{s^k}) \approx \frac k\pi \om_h$ for large $k$. Then 
$\K_2^k\in\ccal^\infty(C\times C)$ for $k$ sufficiently large. 
The diagonal term comes from `self-correlations' of a zero with itself. The second term is the interesting one. In \cite{BSZ00a}, it was shown that $\K_2^k$ has a universal limit using the $1/\sqrtk$ scale of Section~\ref{offdiag}:
\begin{eqnarray}\label{paircor}\tilde K_2^k\left(\frac z\sqrtk,\frac w\sqrtk\right) &\to &\frac{\left(\sinh^2(r^2/2)+r^4/4\right)\cosh(r^2/2) -r^2\sinh(r^2/2)}{\sinh^3(r^2/2)} \\&&= \frac 12 r^2-\frac 1{36} r^6 + \frac 1{720} r^{10}-\cdots,\qquad r=|z-w|\,, \notag\end{eqnarray}
using local holomorphic coordinates about $z_0\in M$ with $\om(z_0)=\frac i2dz\wedge d\bar z$.
Equation \eqref{paircor}, which holds for all Riemann surfaces, was given in \cite{Ha96} for $\CP^1$ and in \cite{NV98} for genus$(C)=1$. 

The fact that the pair correlation $\kappa^k\to 0$ as the distance $r\to 0$ (with $k$ fixed) tells us that the zeros `repel' in the sense that they cluster less than independent random points cluster, as illustrated below:

\begin{center}\includegraphics[height = 2.5
in]{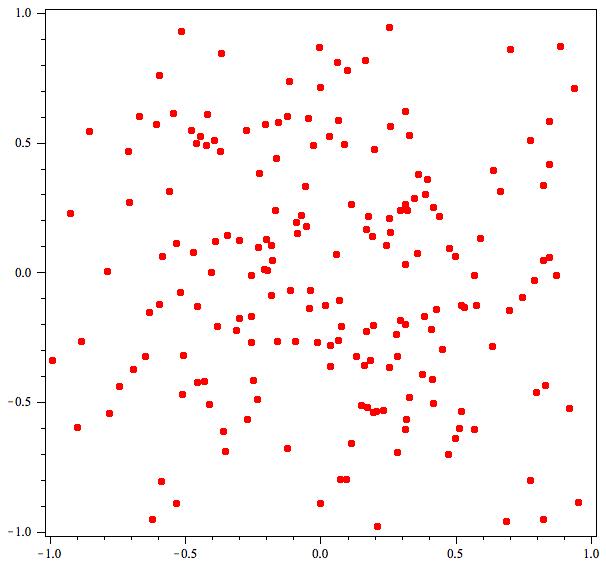}\ 
\includegraphics[height = 2.5in]{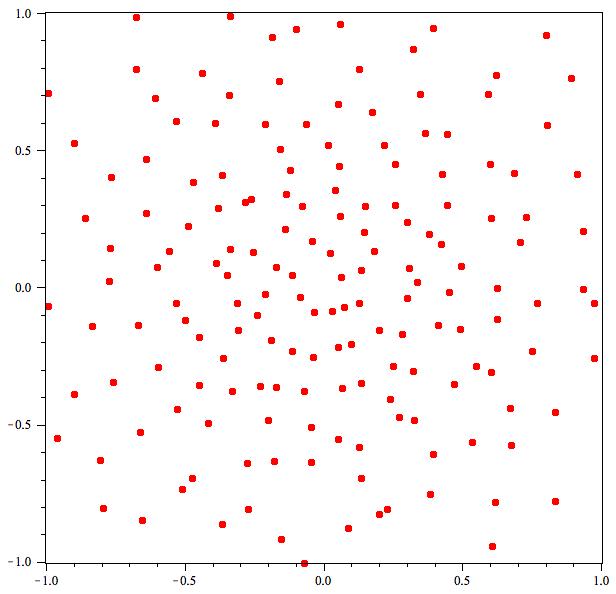}\\{\small\quad\quad
Poisson point process\qquad\qquad\qquad\quad zeros of random polynomials}\end{center}

It was shown in \cite[Th.~3.6]{BSZ00} that $n$-point correlations for random zero sets have universal scaling limits in all dimensions and codimensions of the form
$$\frac 1{k^{np}}\K^k_{npm}\left(\frac{z^1}{\sqrtk}\,,\dots,\,\frac{z^n}{\sqrtk}\right)=\K^\infty_{npm}(z^1,\dots,z^n)+O\left(\frac 1\sqrtk\right)\,, $$ where 
$p$ is the codimension of the simultaneous zero set (of $p$ holomorphic sections of $L^k$). Formulas for $\K^k_{npm}$ are given in \cite{BSZ00} and \cite{BSZ01}. In particular, for the point process case $p=m$,
$$\K^\infty_{2pm}(z,w) = \frac{m+1}{4}\,r^{4-2m}+O(r^{8-2m})\,,\qquad r=\|z-w\|\,.$$
Hence, random simultaneous zeros  of $m$ sections in $H^0(M,L^k)$ do not `repel'  when $\dim M \ge 2$, and in fact for $\dim M\ge 3$ they  cluster more than those in Poisson processes,  for large $k$.

\subsection{A pluri-bipotential for the zero variance}\label{s-bi} \hspace{1cm}
 
   In this section we give a formula for the variance $\var (Z_{s^k})$ of the zero current $Z_{s^k}$ of a Gaussian random holomorphic section $s^k\in H^0(M,L^k)$ (Theorem \ref{BIPOT}).  Let us first describe  the variance of a (general)  random current:

\begin{defin} \label{varcurrent} Let $X:\Om\to \dcal'^j(M)$ be a random variable with values in the space $\dcal_\R'^j(M)$ of real currents of degree $j$ on a manifold $M$.  The variance of $X$ is the current
\begin{equation}\label{vc}
\Var(X): =  \E(X\boxtimes X)-  \E(X) \boxtimes\E(X)\,,\end{equation}
where we use the notation
$$S\boxtimes T = \pi_1^*S \wedge \pi_2^*T \in \dcal'^{p+q}(M\times
M)\;, \qquad \mbox{for }\ S\in \dcal'^p(M),\ T\in \dcal'^q(M)\;.$$
Here, $\pi_1,\pi_2:M\times M\to M$ are the projections to the
first and second factors, respectively. Using more intuitive notation, we shall write
$(S\boxtimes T)(z,w) = S(z) \wedge T(w)$, where $(z,w)$ denotes a point of $M\times M$.\end{defin}
 
The rationale behind Definition \ref{varcurrent} is that the variance of the pairing of $X$ with a compactly supported real test form $\psi\in\dcal_\R^{\dim M - k}(M)$ is given by \beq\label{vv}\var(X,\psi)= \big(\Var(X),\,\psi\boxtimes\psi\big)\,.\eeq

We now show that the variance $\var (Z_{s^k})$ of the zero current depends only on the  normalized \szego kernel $P_k$ given in equation\eqref{Pk}:

\begin{theo} \label{BIPOT} {\rm \cite{ShZ08}} Let $(L,h)\to(M,\om_h)$ be a positive holomorphic line bundle over a compact \kahler manifold.
Then the variance of the zero current of holomorphic sections of $L^k$ is given by
\begin{equation}\label{varcur}{\Var}\big(Z_{s^k}\big)=
-\frac 1{4\pi^2}\d_z\dbar_z\d_w\dbar_w \Li_2[P_k(z,w)^2] \in\dcal'^{2,2}(M\times M)\;,\end{equation}
where $\Li_2$ is the `di-logarithm'
$$ \Li_2(t)= \
\sum_{n=1}^\infty\frac{t^n}{n^2}= -
\int_0^t \frac{\log(1-x)}{x}\,dx\;.$$ 
\end{theo}

Theorem \ref{BIPOT} is equivalent via \eqref{vv}  to the variance formula:
\begin{equation} \label{varint1} \var\left(\int_{ Z_s}\psi\right)=\frac 1{4\pi^2}
 \int_{M\times M}
 \Li_2[P_k(z,w)^2] \;i\ddbar \psi(z) \wedge i\ddbar \psi(w)
\;,\end{equation} for test forms
$\psi\in\dcal^{m-1,m-1}(M)$. (In \cite[Th.~3.1]{ShZ08},   $\wt G(t)=\frac 1{4\pi^2}\Li_2(t^2)$.)

In the Riemann surface case, \eqref{varcur} becomes (in local coordinates)
$${\Var}\big(Z_{s^k}\big) = \frac 1{4\pi^2} \De_z\De_w \Li_2[P_k(z,w)^2]\;(idz\wedge d\bar z)\wedge(idw\wedge d\bar w), $$ so that $\frac 1{4\pi^2}\Li_2[P_k(z,w)^2]$ is a bipotential for the variance of zeros.  In higher dimensions, we say that $ \frac 1{4\pi^2}\Li_2[P_k(z,w)^2]$ is a {\it pluri-bipotential\/} for
the variance current.

To prove Theorem \ref{BIPOT}, we first note that it suffices to
verify the identity over a trivializing neighborhood  $U$ of $(z,w)$.  
{Using the notation of Section \ref{inner+gaussian}, we write $s^k=\langle \vec{c}, \vec{S^k}\rangle$,  where $\vec{S^k}=\vec{F}e_L^{\otimes k}$, $s^k=fe_L^{\otimes k}=\langle\vec c, \vec F\rangle e_L^{\otimes k}$. We have by Theorem \ref{EPL}, \beq\label{EPLf}\E(Z_{s^k})= \frac i{\pi}\ddbar\log\|\vec{F}\|\,.\eeq

The first step of the proof of Theorem \ref{BIPOT} is the following lemma:

\begin{lem} \label{varinta} Writing $\tilde f=\|\vec F\|\inv f$, we have
$$\Var( Z_{s^k}) = -\frac{1}{\pi^2}
\d_z\dbar_z \d_w\dbar_w \E\left(\log |\tilde f(z)|\,\log |\tilde f(w)|\right)\;.$$
\end{lem}

\begin{proof}  By the Poincar\'e--Lelong formula \eqref{PLf},
\beq\label{EZV}\E(   Z_{s^k}\boxtimes  Z_{s^k})= -\frac{1}{\pi^2}\d_z\dbar_z
\d_w\dbar_w \E\left( \log|f(z)|\,\log|f(w)|\right)
\,.\eeq 

Then $\tilde f=\langle\vec{c},\vec u\rangle$, where $\vec u=\|\vec{F}\|\inv\vec F$,  and  
\begin{eqnarray}\log|f(z)|\,\log|f(w)| &=&
\log \|\vec{F}(z)\| \,\log \|\vec{F}(w)\| + \log\|\vec{F}(z)\| \,\log |\tilde f(w)|\nonumber \\&&+ \log|\tilde f(z)| \,\log \|\vec{F}(w)\| + \log|\tilde f(z)| \,\log |\tilde f(w)|\;,\label{4terms}\end{eqnarray} which decomposes (\ref{EZV}) into
four terms. By \eqref{EPLf}, the first term contributes $$-\frac
1{\pi^2}\, \ddbar\log\|\vec{F}(z)\|\wedge
\ddbar\log\|\vec{F}(w)\|=\E( Z_f)\boxtimes\E( Z_f)\,.$$ Since $\|\vec u\|\equiv 1$, $\E(\log |\tilde f(w)|)$ is independent of $w$ and hence
the second term vanishes when applying $\d_w\dbar_w$. The third
term likewise vanishes  when applying $\d_z\dbar_z$. Therefore,
the fourth term yields the variance current $\var( Z_{s^k})$.
\end{proof}

 Next we use the following formula from \cite[Lemma~3.3]{ShZ08}:}

\begin{lem} \label{varintb} Let $(Y_1,Y_2)$ be joint complex Gaussian
random variables of  mean 0 and  variances $\E(|Y_1|^2)= \E(|Y_2|^2)=1$.  Then
$$\E\big(\log |Y_1|\, \log |Y_2|\big) = \frac14\,\Li_2\big(\left|\E(Y_1\overline
Y_2)\right|^2\big)  +\frac{\ga^2}{4}\qquad (\ga = \mbox{Euler's constant})\,.$$\end{lem}

\medskip\noindent{\it Completion of the proof of Theorem \ref{BIPOT}:\/} Fix points $z,w\in
M$, and let  $Y_1=\tilde f(z),\ Y_2=\tilde f(w)$.
Recalling \eqref{2point}, we have
$$|\E(Y_1\overline Y_2)| =  \frac{\sum_jF^k_j(z)\overline{F^k_j(w)}}{\|\vec F(z)\|\,\|\vec F(w)\|} =
\frac{\|B_k(z,w)\|}{\|B_k(z,z)\|^\frac 12\,\| B_k(w,w)\|^\frac
12}= P_k(z,w)\,.$$
Therefore, by Lemma \ref{varintb},
$$\E\left(\log |\tilde f(z)|\,\log |\tilde f(w)|\right) = \E \big(\log |Y_1|\, \log
|Y_2|\big) = \frac14\,\Li_2\big(P_k(z,w)^2\big)  +\frac{\ga^2}{4}\,.$$
Equation \eqref{varcur} then follws from  Lemma \ref{varinta}.
\qed

\subsection{Smooth linear statistics of zeros}
 By  {\it linear statistics} for $H^0(M,L^k)$, we mean the random variable on the probability space $(H^0(M,L),\ga_{h,dV})$
\beq\label{linstat} s^k\ \mapsto \ \left(  [Z_{s^k}], \psi \right) : = \int_{Z_{s^k}} 
\psi (z) ,\qquad s^k\in H^0(M,L^k),\eeq for a fixed continuous test  form $\psi\in\dcal^{m-1,m-1}(M)$. In particular, when $M$ is a Riemann surface $C$, we have 
$$\left(  [Z_{s^k}], f \right)  =  \sum_{z: s^k(z) = 0}f(z),\qquad s^k\in H^0(C,L^k),$$ for a fixed continuous test function $f$.

Both the expectation and the variance of \eqref{linstat} have asymptotic expansions.
To determine the asymptotic expansion of $\E(Z_{s^k}, \psi)$, for $s^k\in H^0(M,L^k)$, we first apply \eqref{seg} to obtain 
\beq \log \|B_k(z,z)\|\sim \log \left(\frac{k^m}{\pi^m}\right) 
+ \frac{\rho_h}{2} k\inv+b_2 k^{-2}+\cdots,\end{equation}
where $\rho_h$ is the scalar curvature of $\om_h$.
Then by Theorem \ref{EPL} and \eqref{seg}, we obtain  the complete asymptotic expansion of the linear statistics  
\begin{equation}\label{ave}\frac 1k\,\E(Z_{s^k},\psi)\sim
\frac 1\pi\int_M \om_h\wedge\psi+\left(\frac i{4\pi}\int_M\rho_h\,\ddbar\psi\right)k^{-2}+\cdots\;. \end{equation}

Similarly, the variance has the following complete asymptotic expansion:
\begin{theo}\label{linearstat}\cite{Sh21}  Let $(L,h)\to(M,\om)$ be a positive holomorphic line bundle over a compact \kahler manifold, and let $\psi\in\dcal_\R^{m-1,m-1}(M)$.  The variance of the linear statistics $(Z_{s^k},\psi)$ has an asymptotic expansion of the form
\beq\var\big(Z_{s^k},\psi\big) \sim A_0k^{-m}+A_1 k^{-m-1} + \cdots +A_jk^{-m-j}+ \cdots.\label{expansion}\eeq The leading and sub-leading coefficients are given by \begin{eqnarray}\label{A0}A_0&=&\frac{\pi^{m-2}\zeta(m+2)}4\,\|\ddbar\psi\|_{2}^2\,,\\A_1&=& 
 -\pi^{m-2}\zeta(m+3)\left\{\frac 18\int_M\rho_h|\ddbar\psi|^2\frac 1{m!}\om_h^m+\frac14\|\d^*\ddbar\psi\|_{2} ^2\right\},\label{A1}\end{eqnarray}  
 where $\zeta$ denotes the Riemann zeta function, and $\rho_h$ is the scalar curvature of $\om_h$. 
 \end{theo}
 The expansion \eqref{expansion} builds on the methods of \cite{ShZ10}, where the asymptotic formula $\var\big(Z_{s^k},\psi\big) = k^{-m}[A_0+O(k^{-1/2+\ep})]$ was given.

In the complex curve case, \eqref{expansion} becomes  \beq\var\left(Z_{s^k}, f\right)\sim \frac {\zeta(3)}{16\pi}\|\Delta f\|^2k\inv -\frac{\pi^3}{2880}\left\{\int_M\rho_h|\Delta f|^2\,\omega+\|d\Delta f\|_{2} ^2\right\}k^{-2}+\cdots,\label{sharpRiemann}\eeq
for $f\in\ccal^\infty(M)$.

Thus, smooth linear statistics  are self-averaging in the sense  that its
fluctuations  are of smaller order than its typical values. The fact that
the variance involves $\|\Delta f\|^2_2$ rather than $\|\nabla f\|^2_2$
signals that the covariance kernel is not $\Delta^{-1}$ but $\Delta^{-2}$.

\subsection{Asymptotic normality of zero  distributions}\label{s-normality}

The following theorem was proved first by Sodin--Tsirelson \cite{ST04} for certain model
random analytic functions on $\C, \CP^1$ and the unit disc and then  in 
\cite{ShZ10}   to general one-dimensional ensembles and to
codimension one zero sets in higher dimensions:
\begin{theo} \label{AN} \cite{ShZ10} Let $(L,h)\to(M,\om)$ be a positive holomorphic line bundle over a compact \kahler manifold, and let $\psi$ be a real  $(m-1,m-1)$-form on $M$ with
$\ccal^3$ coefficients.  Then for random sections $s^k\in H^0(M,L^k)$, the distributions of the random variables
$$
k^{m/2}\left(Z_{s^k}-\frac k\pi\,\om,\psi\right)$$ converge weakly to the Gaussian distribution of mean 0 and variance
$\frac{\pi^{m-2}\,\zeta(m+2)}{4}\, \|\ddbar\psi\|_2^2$, as $ k \to \infty$.\end{theo}

Theorem \ref{AN} follows from a  general result of Sodin-Tsirelson \cite{ST04} on asymptotic normality
of nonlinear functionals of Gaussian processes and the properties of the normalized \szego kernel \eqref{Pk}.
To describe the result of \cite{ST04}, we recall that a (complex) Gaussian process on a measure space $(T,
\mu)$ is a random variable (with values in the space of complex measurable functions on $T$) of the form
$$w(t)=\sum c_j  g_j(t)\;,$$ where the $c_j$ are i.i.d.\ complex
Gaussian random variables of mean 0, variance 1, and the $g_j$ are
(fixed) complex-valued measurable functions. We say that $w(t)$ is {\it normalized} if $\sum|g_j(t)|^2=1$ for all $t\in T$; i.e., if $w(t)\sim \ncal_\C(0,1)$ for all $t\in T$.

\begin{theo}\label{ST} \cite{ST04} 
Let $w^1,w^2,w^3,\dots$ be a sequence of normalized complex
Gaussian processes on a finite measure space $(T, \mu)$. Let $f:\R^+\to\R$ be monotonically increasing such that $f(r)
\in L^2(\R^+, e^{-r^2/2} rdr)$,  and let
$\eta:  T
\to
\R$ be bounded measurable.

Let $\ccal_k(s, t):= \E\left(w^k(s)\overline{w^k(t)}\right)$ be the
covariance function for $w^k$ and suppose that
\begin{enumerate}

\item[i)]  \quad $\displaystyle\liminf_{k \to \infty} \frac{\int_T
\int_{T} |\ccal_k(s, t)|^{2 } \eta(s) \eta(t) d\mu(s)
d\mu(t)}{\sup_{ s \in T} \int_T |\ccal_k(s, t)| d\mu(t)} > 0\;;$\
\item[ii)]  \quad $\displaystyle\lim_{k \to \infty}\; \sup_{s \in T}
\int_T |\ccal_k(s, t)| d\mu(t) = 0.$

\end{enumerate}

 Consider the random variables
$$Y_k = \int_T f(|w^k(t)|) \eta(t) d\mu(t).$$

Then the distributions of the random variables
$$ \frac{Y_k - \E Y_k}{\sqrt{\var(Y_k)}}$$
converge weakly to
$\ncal(0, 1)$ as $ k \to \infty$.
\end{theo}

To prove Theorem \ref{AN}, we apply Theorem \ref{ST} with $f(r)=\log r$ and $(T, \mu) = (M,
dV)$. To define the normalized Gaussian processes $w^k$ on $M$, 
choose a measurable section $\sigma_L:M\to L$ of $L$ with
$\|\sigma_L(z)\|_h=1$ for all $z\in M$, and  let
$$S_j^k=F^k_j\sigma_L^{\otimes k}\,,\quad j=1,\dots,d_k,$$ be an
orthonormal basis for $H^0(M, L^k)$.  We then  let
$$w^k(z) : =\sum_{j=1}^{d_k}c_j\, \frac{F^k_j(z)}{\sqrt{\Pi_k(z,z)}}\,.$$
Since $|F^k_j|=\|S^k_j\|_{h^k}$, it follows that $w^k$
defines a normalized complex
Gaussian process. In fact,
\beq\label{ncgp}\sqrt{\Pi_k(z,z)}\,w^k\sigma_L^{\otimes k}= \sum c_j S^k_j = s^k\,,\eeq
where $s^k$ is a random holomorphic section in $ H^0(M, L^k)$. 

We now let $\psi$ be a fixed   real
$(m-1,m-1)$-form on $M$ and we write
$$\frac i\pi\ddbar\psi =\eta\,dV\;.$$ 
Then by \eqref{ncgp},
$$Y_k= \int_M\log |w^k| \eta\,dV= \int_M\left( \log \|s^k(z)\|_{h^k} - \log
\sqrt{\Pi_k(z,z)} \right)\frac i\pi \ddbar \psi(z) = \big(
Z_{s^k},\psi\big) + C_k\,,$$
where each $C_k$ is a constant
independent of the random section $s^k$. Hence $Y_k$ has the
same variance  as the linear statistic $\big(Z_{s^k},\psi\big)$.
In fact, the covariance functions
$\ccal_k(z,w)$ for these Gaussian processes satisfy
$$|\ccal_k(z,w)|=P_k(z,w)\;.$$

It was shown in \cite{ShZ10}, using the properties of the normalized Bergman kernel $P_k$ given by Propositions~\eqref{DPdecay}--\eqref{better}, that conditions (i)--(ii) of Theorem \ref{ST} hold. Hence, the distributions of the random variables
$$\frac{(Z_{s^k},\psi)-\E(Z_{s^k},\psi)}{\sqrt{\var(Z_{s^k},\psi)}}= \frac{Y_k-\E(Y_k)}{\sqrt{\var(Y_k)}}$$
converge weakly to the
standard Gaussian distribution
$\ncal(0, 1)$ as $ k \to \infty$.

The conclusion of Theorem \ref{AN} then follows from the leading asymptotics of the expectation $\E(Z_{s^k},\psi)$ and the variance $\var(Z_{s^k},\psi)$ given by equations \eqref{ave} and \eqref{expansion}.\qed

Nazarov and Sodin \cite{NS11, NS12} give results on variances and asymptotic normality  for linear statistics on $\C$ with test functions that are not continuous. It is open whether similar results hold for line bundles on compact \kahler manifolds and whether asymptotic normality holds for linear statistics of zeros of any codimension.

One can also consider linear statistics for the point process consisting of simultaneous zeros of $m$ independent random sections of $H^0(M,L^k)$:
\beq\label{linstatm}(s^k_1,\dots,s^k_m)\mapsto \left(  [Z_{s^k_1,\dots,s^k_m}], f \right)  =  \sum_{\{ s_1^k(z) =\cdots=s_m^k(z)= 0\}}f(z),\eeq for a fixed continuous test function $f$.
 Asymptotics for the expectation and variance of \eqref{linstatm} in \cite{ShZ10}, but it is an open problem with asymptotic normality holds for  \eqref{linstatm}.

\subsection{Counting random zeros in a set and hole probabilities}

We say that  a set $U\subset M$ is a `hole' in the zero set $Z_s$ of a section $s\in H^0(M,L)$ if  $  Z_{s}\cap U= \emptyset$. Hole probabilities,
overcrowding and other number statistics for special ensembles of functions of one
complex variable were given, for example, in \cite{Kr,ST2}. To describe the
framework of  these results in the case where $\dim M=1$, i.e. where $M$
is a compact Riemann surface $C$, we consider the random variable $\ncal^U_k(s^k): = \#(Z_{s^k}\cap U)$
on $H^0(C, L^k)$  which counts the zeros  of a section $s^k$ in an open set $U$. Hence the hole probability is the probability that
$\ncal^U_k(s^k) = 0$.  Sodin and Tsirelson \cite{ST2} gave an asymptotic formula 
for the variance of $\ncal^U_k$ when
$M=\CP^1$ (and for the analogous cases of holomorphic functions on $\C$ and on the
disk).  This formula was  sharpened and generalized using Theorem~\ref{BIPOT} to
arbitrary compact Riemann surfaces  as well as to compact \kahler manifolds of
any dimension, to obtain the following analogue of Theorem~\ref{linearstat}:
\begin{theo}\label{numbervar} \cite{ShZ08} Let $(L,h)\to(M,\om_h)$ be a positive holomorphic line bundle over a compact \kahler manifold. Let $U$ be a domain in $M$ with piecewise $\ccal^2$ boundary without cusps. Then for $m=\dim M$ independent Gaussian random sections $s^k_1,\dots,s^k_m$  in $H^0(M,L^k)$, the variance of the random variable 
$$\ncal^U_k(s^k_1,\dots,s^k_m): = \#\{z\in U:s^k_1(z)=\cdots=s^k_m(z)=0\}$$ has the asymptotics
\beq\label{zerovar} \var(\ncal^U_k)=k^{m-1/2}[\nu_{m}\vol_{2m-1}(\d U)+O(k^{-1/2+\ep})]\,,\eeq where $\nu_{m}$ is a universal positive constant depending only on $m$.\end{theo}
For the Riemann surface case, $\nu_{1}=\zeta(3/2)/(8\pi^{3/2})$.
For random zero sets of one section, we similarly have:
\begin{theo}\label{numbervar2} \cite{ShZ08} With the hypotheses of Theorem \ref{numbervar},
$$\var(\vol_{2m-2}Z_{s^k}\cap U)= k^{3/2-m}\left[{\textstyle\frac18}\pi^{m-5/2}\zeta(m+1/2)\,\vol_{2m-1}(\d U)+O(k^{-1/2+\ep})\right].$$\end{theo}

Theorems \ref{numbervar} and \ref{numbervar2} are special cases of a general result for simultaneous zeros of $p$ holomorphic sections on $M$, for $1\le p\le m$.

The volume of the zero set $Z_{s^k}$ inside a domain also satisfies  a large deviations bound of the form:

\begin{theo}\label{volumes} \cite{ShZZ08} Let $(L,h)\to (M,\om_h)$ be as in Theorem
\ref{numbervar}, and let  $U$ be an open subset of $M$ such that $\d U$ has zero measure in $M$.
Then for all
$\de>0$ sufficiently small, there is a constant $C_{U,\de}>0$  such
that
$$Prob\left\{ \left| \vol_{2m-2}(Z_{s^k}\cap U) - \frac m\pi\,
\vol_{2m}(U)\,k
\right| > 
\delta \,k\right\}\leq e^{- C_{U,\delta} k^{m+1}}\quad \forall\ k\gg 0\,.$$
\end{theo}

Here, $k\gg 0$ means $k\ge k_0$ for some $k_0\in \Z^+$. 
In particular, for the case where $\dim M=1$, we have:

\begin{cor}\label{dim1} 

 Let $(L,h)\to (C,\om_h)$ be a positive Hermitian line bundle over a compact Riemann
surface.  Let
$U\subset C$ be an open set in $C$ such that $\d U$ has zero measure in
$C$.  Then for all $\de>0$, there is a constant $c_{U,\de}>0$ such that 
$$ Prob \left\{ s^k:\;\left| \ncal^U_k(s^k) - \frac k\pi \,\mbox{\rm Area}(U)
\right|> 
\delta\,k\right\} \leq e^{- c_{U, \de} k^2}\quad \forall\ k\gg 0\,.$$
\end{cor}

We also have upper and lower estimates for the `hole
probability':

\begin{theo}\cite{ShZZ08}\label{hole} Let $(L,h)\to (M,\om_h)$ and $U\subset M$ be as above, and suppose there is a section $s\in H^0(M,L)$ that does not vanish anywhere on $\overline U$. Then there exist  constants $c_U'>c_U>0$ such that 
\beq\label{twosided}e^{- c'_{U} k^{m+1} } \leq Prob
\{s^k:  Z_{s^k}\cap U=\emptyset\}
\leq e^{- c_{U} k^{m+1} }\quad \forall\ k\gg 0\,. \eeq \end{theo}
The upper bound in \eqref{twosided} is an immediate consequence of Theorem \ref{volumes} with $\de<\frac m\pi
\vol_{2m}(U)$. The analogue of \eqref{twosided} for random entire functions of one variable was given in \cite{ST2}.

It is an open question whether $$\log(Prob
\{s^k:  Z_{s^k}\cap U=\emptyset\}) \sim {-\tilde c(U) k^{m+1}}$$ for some constant $\tilde c(U)$. This was shown in \cite{Zhu14}  to hold for $M=\CP^m$ and $U=\Delta_r^m$, where   $\Delta_r=\{z\in\C:|z|<r\}$, with a  specific formula for $\tilde  c(\Delta_r^m)$. When $r\ge 1$, $$\tilde  c(\Delta_r^m)=\frac{2m\log r}{(m+1)!}+\frac 1{m!}\sum_{j=1}^m\frac 1{j+1}\,.$$

\subsection{Expected local topology of random zero sets} 

An active topic of recent research in random {\it real} algebraic geometry is the random topology of random real algebraic varieties defined by the zero
locus   of one or several independent real polynomials of a fixed degree $k$: the number of connected components, the betti numbers, and the combinatorics
of the components. Works of  Gayet and Welschinger \cite{GaW14,GaW15,GaW17} and others have resolved many problems in this area. For this survey,
the question is whether there exist problems of this nature in the complex case. Globally, the answer is no: all of these topological invariants are deterministic. 
Recently,  D. Gayet \cite{G22} has studied the local analogues of these problems in open sets $U \subset M$, and the local  topology of
$Z_{s^k} \cap U$ is very random. The main result of \cite{G22}, stated here in the rank 1 case, is the following:
\begin{theo}\label{GAYETTh1} Let $(L,h)\to(M,\om_h)$ be a positive holomorphic line bundle over a compact \kahler manifold. Let $U$ be  an open set  in $M$ with smooth boundary. Then for random holomorphic sections $s^k\in H^0(M,L^k)$,
$$  \lim_{k \to \infty}  \frac{1}{k^m} \E\, b_j ( Z_{s^k} \cap U) = \left\{ \begin{array}{cl} 0, & \mbox{\rm for }\  0\le j\le 2m-2,\ j\neq m-1 \\ [6pt]\frac {m!}{\pi^m}\,\vol_{2m}(U)\,,& \mbox{\rm for }\ j=m-1\end{array}\right..
$$
\end{theo} 

By contrast, in the real domain, all Betti numbers grow like $\lambda^m$ where $1/\lambda$ is the natural scale of the model. In the complex setting, $\lambda = \frac{1}{\sqrt{k}}$. 
Yet the $(m-1)$-th Betti number (and only that Betti number) grows like $k^m$.

In \cite{Au97}, Auroux  proved that the  (deterministic) quantitatively transversal Donaldson hypersurfaces, which are zeros of sections that vanish transversally with a controlled derivative, satisfy this local topology estimate for the $(m-1)$-th Betti number.  The Donaldson hypersurfaces have common features with the random ones of Theorem \ref{GAYETTh1} in the complex setting. For instance, the current of integration  $Z_s$ fills out $M$  uniformly for large degrees $k$ in both contexts.

Gayet explains the intuition behind Theorem \ref{GAYETTh1} when $m=2$,  as follows: By the maximum principle, if a complex curve  $C \subset \C^2$ contacts a real hyperplane $H$, and is locally on one side of $H$, then $C \subset H$. Let $u: \C^2 \to \R$ be a Morse function and for $r> 0$ let
$u_r(z) = u(r^{-1}z)$. For increasing  $r$, the level sets of $u_r$ locally become closer and closer to being planar so that there are fewer  random curves touching them from 
the interior; i.e., there are fewer  critical points of $u |{Z_{s^k}}$ of index $0$ compared to critical points of index $1$. The result then follows from Morse theory.

As this intuition suggests, the proof of Theorem \ref{GAYETTh1} involves the strong Morse inequalities and  a statistical study of critical points via the Kac--Rice formula.

\section{Critical points and  values of random holomorphic sections}
 In this section, we review some results on random critical points, critical values, and excursion sets of Gaussian random holomorphic sections of $H^0(M, L^k)$ and their asymptotics. 

\subsection{Critical points}\label{critpoints}

 The   critical point
set $\crit(s, h)$ of a holomorphic section $s\in H^0(M,L)$ is defined by
\begin{equation}  \label{Critdef} \crit(s, h): = \{z\in M: \left.d (\|s\|_h^2)\right|_z = 0,\; s(z) \not= 0\}= \{z\in M: (\nabla_h s)(z) =
0 ,\; s(z) \not= 0\},
\end{equation}
where $\nabla_h$ is the Chern connection of the Hermitian holomorphic line bundle $(L,h)$. 

Along with the expected distribution of zeros (of a random section of a line bundle on a Riemann surface, or of simultaneous zeros in higher dimensions), we are also interested in the expected distributions   of  critical points 
\begin{equation}  \K^{k\,\crit}_{1} :=\textstyle \E \left[\sum_{z\in \crit (s^k, h^k)} \delta_{z}\right]
\label{Kcrit}\end{equation}
of  random holomorphic sections $s^k\in H^0(M,L^k)$. 

Additionally, the total number of critical points, $\#\crit(s,h)$,  is a (nonconstant) random variable, unlike the total number of zeros of $m$ holomorphic sections of $L\to M$, which equals the topological invariant $c_1(L)^m$. Although the alternating sum of critical points of each Morse index is a topological invariant,  the sum $\#\crit(s,h)$ as well as the number of critical points of a given Morse index is non-constant. Hence, we are interested in the average number of critical points,
\begin{equation}\ncal^{\crit}_{k}(L,h):=\E\left(\#\crit(s^k,h^k)\right)   =\int_M
d\K^{k\,\crit}_{1}\,, \end{equation} of a random section $s^k\in H^0(M,L^k)$. We also consider the average number of critical points of Morse index q, which we denote by $\ncal^{\crit}_{k,q}(L,h)$, for $m\le q\le 2m$.  (It was observed in  \cite{DSZ04} that the critical points are all of Morse index $\ge m$.) 

The expected number $\ncal^{\crit}_{k}(L,h)$ and expected distribution of critical points can be expressed as  formulas involving the Bergman kernel (see \cite[Theorems~1\,\&\,6]{DSZ04}).  These formulas can be used to  obtain explicit expressions for the hyperplane section bundle $\ocal(1)\to\CP^m$ with the Fubini-Study metric. In fact, the averages $\ncal^\crit_k(\ocal(1)\to\CP^m)$ are rational functions of $k$. In dimensions 1 and 2, we have 
\begin{eqnarray*}{\ncal}^\crit_{k}\left(\ocal(1)\!\to\!\CP^1\right)&=& \frac{5k^2-8k+4}{3k-2}\;,\\
\ncal^\crit_k\left(\ocal(1)\!\to\!\CP^2\right)&=& 
{\frac {59\,{k}^{5}-231\,{k}^{4}+375\,{k}^{3}-310\,{k}^{2}+132\,k-24}{
 \left( 3\,k-2 \right) ^{3}}}\;.\end{eqnarray*}
The average numbers of critical points of sections of $\ocal(k)\to\CP^m$ of each Morse index $q$ ($m\le q\le 2m$) are likewise rational functions of $k$, for  $m\ge 1$.  See \cite[Appendix~1]{DSZ06} for some explicit formulas.

The expected number of critical points depends on the metric $h$, but surprisingly its leading asymptotics is independent of $h$.  This can be seen from the following asymptotic formulas for $\ncal^{\crit}_{k}(L,h)$ and $\ncal^{\crit}_{k,q}(L,h)$: 

\begin{theo} \label{critk} {\rm (\cite[Cor.~1.4]{DSZ06})} Let $(L,h)\to (M,\om_h)$ be a
positive holomorphic line bundle on a compact \kahler manifold. Then for $m\le q\le 2m$,\begin{multline}\label{crit-exp}\ncal^{\crit}_{k,q}(L,h) \sim
\left[\be_{0qm}\,c_1(L)^m\right]k^m +
\left[\be_{1qm}\, c_1(M)\cdot
c_1(L)^{m-1}\right]k^{m-1} 
\\+ \left[\be_{2qm}\int_M\rho_h^2\,dV_h + \be'_{2qm}\, c_1(M)^2\cdot
c_1(L)^{m-2} + \be''_{2qm}\, c_2(M)\cdot c_1(L)^{m-2}\right]
k^{m-2}+\cdots\,,\end{multline} where $\rho_h$ is the scalar curvature of $\om_h$, and
$\be_{0qm},\be_{1qm},\be_{2qm},\be'_{2qm},\be''_{2qm}$ are universal constants
depending only on $q$ and $m$.
\end{theo}

Taking the sum over $q$ from $m$ to $2m$, one obtains a similar asymptotic expansion for $\ncal^{\crit}_{k}(L,h)$.

In fact, both the leading term and  the subleading term of the expansion do not depend on the choice of metric on $L$. Note that the   non-topological part of the third term is $\be_{2qm}$ times the Calabi functional $$Cal(\om_h):=\int_M\rho_h^2\,dV_h\,.$$ 
It was shown in \cite{Bau10} that $\sum_{q=m}^{2m}\be_{2qm}$ is positive for all $m\ge 1$, and hence we have:
\begin{cor} \label{a-minimizing} Let $h$ and $h'$ be metrics on $L$ with positive curvature. If $Cal(h')>Cal(h)$, then there exists $k_0$ such that $\ncal^{\crit}_{k}(L,h')>\ncal^{\crit}_{k}(L,h)$ for all $k\ge k_0$.\end{cor}

It is known that  \kahler metrics of constant scalar curvature are  critical metrics of  the functional $Cal$ on the space of \kahler metrics in $c_1(L)$ and that all critical metrics are global minimums for $Cal$ on this space \cite{Ca85, Hw95}. Furthermore, if $c_1(L)$ contains a constant scalar curvature \kahler metric, then every critical metric has constant scalar curvature \cite{Ca85}. Thus, constant scalar curvature metrics have ``asymptotically minimal critical numbers" in the sense that if $h$ and $h'$ are two positive metrics on $L$  such that $\om_h$ but not $\om_{h'}$ has constant scalar curvature, then the conclusion of Corollary \ref{a-minimizing} holds.

It is known that $\be_{2mm}>0$ for all $m$ \cite{Bau10} and that  $\be_{2qm}>0$ in low dimensions \cite{DSZ06}, but it remains  an open question whether $\be_{2qm}>0$ for all $q,m$.

Like the  zeros of a random section on a Riemann surface (or the simultaneous zeros of $m$-sections on $M$), the critical points of random sections on $M$ form a point process.  Its one-point function (or expected density) has an asymptotic expansion \cite{DSZ06}
\beq\label{crit-ave} \frac1{k^m}\K^{k\,\crit}_1(z) \sim \left[b_0+b_1(z)k^{-1}+b_2(z)k^{-2}+\cdots\right]\frac1{m!}\om_h^m\,.\eeq  The constant $b_0$ depends only on the dimension $m$; see \cite{DSZ06} for its values. To our knowledge, the only result on the pair correlation of critical points of sections in $H^0(M,L^k)$ is  the asymptotic formula for the case of Riemann surfaces in \cite{Bab12}:
\begin{theo}Let $(L,h)\to(C,\om_h)$ be a positive line bundle over a compact Riemann surface $C$. Then the pair correlation of critical points $\tilde K^{k\,\crit}_2(z,w)$ has the scaling limit
$$\lim_{k\to\infty}\frac 1{k^2}\, \tilde K^{k\,\crit}_2\left(\frac z\sqrtk,\frac w\sqrtk\right) = \frac 2{3\pi^2}+O(r^2)\,,$$ where $r=\dist(z,w)$. \end{theo}
Thus the clustering of critical points on Riemann surfaces is similar to that of Poisson point processes.  It is an open problem to determine critical point pair correlation formulas in complex dimension greater than one.

\subsection{\label{EXCURSSECR}Sup norms and random excursion sets} In this section, we discuss the excursion sets$$\{z\in M:\|s(z)\|_h>r\} $$  and sup norms of  random holomorphic sections $s$  of  positive holomorphic line bundles  $(L, h)$. Random excursion sets have been much studied in the case of real Gaussian fields for over three decades; we refer to \cite{G21, TA03}  for background and references to that subject. 
 
It is most useful to study the sup norms and excursion sets of random sections of unit $\lcal^2$ norm with respect to  Haar probability measure $\nu_{k}$ on the unit $(2d_k-1)$-sphere, \begin{equation} SH^0(M, L^k) =\{s^k \in H^0(M, L^k) : \|s^k\|_{G(h^k,dV)} = 1\},
\end{equation} which probability space we call the {\it spherical ensemble}. 
One well-known problem is to determine the expected  Euler characteristics of the excursion sets,  $$\E_{\nu_k}\big[\chi\{\|s^k(z)\|_{h^k}>r\}\big]=\int_{H^0(M,L^k)}\chi\{\|s^k(z)\|_{h^k}>r\}\,d\nu_k(s), $$ and also the probability that the excursion set is non-empty. (Here and in the following, $\chi$ denotes the Euler characteristic.)

We  note that  for all $s^k=\sum c_jS^k_j\in SH^0(M,L^k)$, we have by \eqref{bergmankernel}--\eqref{seg},
\beq\label{upperbound} |s^k(z)|^2_{h^k}\le \sum_{l=1}^{d_k}|c_l|^2\sum_{j=1}^{d_k} \|S^k_j(z)\|^2_{h^k}=1\cdot \|B_k (z,z)\| =\left(\frac1{\pi^m}+o(1)\right)k^m,\eeq
and thus the excursion sets $\{z:\|s^k(z)\|_{h^k}> (\pi^{m/2}+\ep)k^{m/2}\}$ are  empty for all $s^k\in H^0(M, L^k)$, for $k$ sufficiently large.
Furthermore, a much sharper asymptotic upper bound usually holds: 
\begin{theo} \cite[Th.~1.1]{ShZ03b} For all $n>0$, there exists a positive constant  $C_n$  such that the probability 
$$Prob\left\{\sup_{z\in M}\|s^k(z)\|_{h^k}>C_n\,\sqrt{\log k}\right\}<O\left(\frac{1}{k^n}\right),\qquad s^k\in H^0(M,L^k).$$\end{theo}
Thus for almost all random sequences $\{s^k\}\in\prod_{k=1}^\infty SH^0(M, L^k)$, \beq\label{asbound}\|s^k\|_\infty =O\left(\sqrt{\log k}\right).\eeq

 To our knowledge, the only article studying excursion sets in the holomorphic setting is  a paper of Jingzhou Sun \cite{Sun12}, and we summarize Sun's results below. 

\begin{theo}\label{Euler}
Let $(L,h)\to (M,\om_h)$ be a
positive holomorphic line bundle on a compact \kahler manifold, and let $k_0$ be sufficiently large so that $L^k$ is very ample for all $k\ge k_0$.
Then there exists $u_0<1$ independent of $k$, such that for  for $1\ge u > u_0$,  and $k\ge k_0$,
\begin{enumerate} \item[i)] the excursion set $$\ecal^k_u(s^k):=\left\{z\in M:\|s^k(z)\|_{h^k}>u\,\|B_k(z,z)\|^{1/2}\right\}$$ is either empty or contractible, and \vspace{-26pt}
\item[ii)] $\begin{array}{lcl}\\[18pt] Prob\big[\ecal_u^k(s^k)\neq\emptyset\big]&=&\displaystyle  Prob\left\{\sup_{z\in M}\|s^k(z)\|_{h^k}>u\,\|B_k(z,z)\|^{1/2}\right\}\\[12pt]
&=&\displaystyle \int _M
c(M)(1-kc_1(L))\wedge(ku^2\,c_1(L)-u^2+ 1)^{d_k-1}\,,
\end{array}$\\[8pt]
where  $c(M)(1-kc_1(L))$ is the Chern
polynomial evaluated at $1-kc_1(L)$.
\end{enumerate}\end{theo}
\smallskip

\begin{cor}\label{Euler-est}
With the hypotheses and notation of Theorem \ref{Euler}, for $1\ge u > u_0$, the expected Euler characteristic
\begin{eqnarray*}\E\,\chi\big[\ecal_u^k(s^k)\big] &=& \int _M
c(M)(1-kc_1(L))\wedge(ku^2\,c_1(L)-u^2+ 1)^{d_k-1}\\ &=& (1+o(1))d_k^{m+1}(1-u^2)^{d_k-m-1}u^{2m}\,.
\end{eqnarray*}
\end{cor}
Recall that $$d_k=\dim H^0(M,L^k)=\frac{k^m}{m!}{\int_Mc_1^m(L)}+O(k^{m-1})\,,$$
by the Hirzebruch--Riemann--Roch formula and Kodaira vanishing theorem (or  see e.g., \cite[Lemma~7.6]{ShS85}).

When $\dim M=1$,  Theorem \ref{Euler} yields the following:

\begin{cor}\label{Euler1}
Let $(L,h)$ be a positive line bundle of degree $\de$ over a compact Riemann surface $(M,\om_h)$ of genus $g$.   Then with the notation of Theorem \ref{Euler}, there exists $u_0<1$ independent of $k$, such that for $1\ge u > u_0$, the expected
Euler characteristic
$$ \E\,\chi\big[\ecal_u^k(s^k)\big] =(1-u^2)^{k\de - g-1}\,\big[k^2\de^2u^2-k\de (gu^2-1+u^2)+(2-2g)(1-u^2)\big]
\,,$$
for $\,k\de>2g-2$.
\end{cor}

To  prove Theorem \ref{Euler}, J. Sun proves an embedding theorem of independent interest:

\begin{theo}\label{critical} Let $(L,h)\to (M,\om_h)$ be a
positive holomorphic line bundle on a compact \kahler manifold, and let $k_0\in\Z^+$ so that $L^k$ is very ample for all $k\ge k_0$.
Let $\Phi_k:M\rightarrow \CP^{d_k-1}$ be an embedding given by an orthonormal basis of $H^0(M,L^k)$ with respect to the Hermitian inner product $G(h^k,dV)$, for $k\ge k_0$. Let $r_k$ be the critical radius of $\Phi_k(M)\subset \CP^{d_k-1}$, where $\CP^{d_k-1}$ is given the Fubini--Study metric. Then $\inf_{k\ge k_0} r_k>0$.
\end{theo}
\smallskip

The proof of Theorem \ref{Euler} uses the volume of tubes formula in \cite{G85}, setting $u_0= \cos r_0$, where $r_0=\inf_{k\ge k_0} r_k$. 

\subsection{Critical values}\label{critvalues}

We now turn to the distribution of critical values of random holomorphic sections of powers of $L\to M$. By the ``value" of $s^k(z)\in H^0(M,L^k)$, we mean the norm $\|s^k(z)\|_{h^k}\in \R^+$. We study the norms, since the values $s^k(z)$ lie in different fibers  $L^k_z$ of $L^k$. Thus we let
\begin{equation}  {\rm CV}({s^k}) := \{ \|s^k(z)\|_{h^k} : z\in\crit(s^k,h^k) \} = \{\|s^k(z)\|_{h^k} :z\in M\,,\nabla s^k(z) = 0\,, s^k(z)\neq 0\}
\label{Nabla} \end{equation} denote the set of  critical values of a section $s^k\in H^0(M,L^k)$.

Since  ${\rm CV}(\la{s^k})=|\la|\,{\rm CV}({s^k})$, it is most useful to study the distribution of critical values of random sections of unit $\lcal^2$ norm with respect to  Haar probability measure $\nu_{k}$ on the unit $(2d_k-1)$-sphere \begin{equation} SH^0(M, L^k) =\{s^k \in H^0(M, L^k) : \|s^k\|_{\lcal^2} = 1\},
\end{equation} which probability space we call the {\it spherical ensemble}. (Clearly, the expected distribution of critical points of $s^k$ in the spherical ensemble is identical to that in the Gaussian ensemble and both ensembles have the same expected numbers $\ncal^{\crit}_{k,q}(L,h)$.)

We recall that by \eqref{upperbound}, we have the deterministic bound ${\rm CV}(s^k) \subset (0,Ck^{m/2})\,,$ and furthermore by \eqref{asbound}, ${\rm CV}(s^k) \subset\left(0,C\sqrt{\log k}\right)$ almost surely.

We let 
\beq\label{critvalmeasure} \textstyle [{\rm CV}({s^k})]= \sum_{z\in\crit(s^k,h^k)}\de_{|s^k(z)|_{h^k}}\eeq
denote the critical value distribution of a section $s^k\in SH^0(M,L^k)$.
To describe the asymptotics of the spherical averages  $\E_{\nu_k}[{\rm CV}({s^k})]$, we use the following notation: denote by $\sym(m,\C)\cong \C^\frac{m^2+m}{2}$ the space  of $m\times m$ complex symmetric
matrices, and  define the special (positive definite)  operator   \begin{equation}
\label{Qmatrix} Q=(Q_{jq}^{j'q'}):=\big(\delta_{jj'}\delta_{qq'}+\delta_{jq'}\delta_{qj'}\big)\,,\quad 1\le j\le q \le m, \ 1\le j'\le q'\le m. \end{equation}
We  define the universal function  (depending only on the dimension $m$)
$$f_m(t)=\frac {2}{\pi^{m(m+3)} } \int_{\sym(m,\C)} e^{-|\xi|^2}\left|\det\left(\left|\sqrt Q \,\Xi\right|^2-t^2I\right)\right|d\Xi, $$ where $d\Xi$ denotes Lebesgue measure on $\sym(m,\C)\cong \R^{m^2+m}$.   We then have:

\begin{theo} \label{sphericallimit} \ {\rm (\cite{FZ14a}) } Let $(L,h)\to (M,\om_h)$ be a
positive holomorphic line bundle on a compact \kahler manifold. The normalized expected density of critical values in the spherical ensemble $SH^0(M, L^k)$  has the asymptotics $$ \frac 1{k^m}\E_{\nu_k}[{\rm CV}({s^k})] \to \vol(M)\,f_m(t)\,te^{-t^2}\,dt\,. $$ \end{theo}
In the case of complex curves,  $f_1(t)=\frac 1\pi(2t^2-4+8e^{-{t^2}/2})$.

\section{Point processes and  \kahler metrics}\label{pp}

In this section we give two examples showing how  \kahler metrics can be constructed using point processes.

\subsection{Zero point processes}\label{ppm}

Given a positive line bundle $(L,h)$ over a compact Riemann surface $C$, we can form the point processes $Z_{s^k}$ of zeros of random sections of powers $L^k$ of the line bundle. Recall from Theorem \ref{EZsk} that the expected measure $\frac 1k \E Z_{s^k}$ converges to $\frac 1\pi \om_h$. In higher dimensions, for holomorphic sections ${s_1,\cdots,s_m}$ of a line bundle $L\to M$, chosen so that their common zero set 
$$Z_{s_1,\dots,s_m}=\{z\in M:s_1(z)=\cdots=s_m(z)=0\}=\{\zeta_1,\cdots,\zeta_p\}$$ is finite, we define the empirical probability measure $$\frac 1p[Z_{s_1,\dots,s_m}]=\frac 1p\sum_{j=1}^p \delta_{\zeta_j}\,.$$
Given random holomorphic sections of a positive line bundle $L^k\to M$, the probability measures $\gamma_{h^k,dV}$ on $H^0(M,L^k)$ induce  point processes $Z_{s^k_1,\dots,s^k_m}$ on $M$ (for $k$ sufficiently large so the common zero set is 0-dimensional a.s.). We then have

\begin{theo} \label{volume-limit}\cite[Prop.~4.4]{ShZ99} Let $(L,h)\to (M,\om_h)$ be a positive holomorphic line bundle over a compact \kahler manifold. Then 
$$\frac 1{k^m}\,\E [Z_{s^k_1,\dots,s^k_m}] \to \frac 1{\pi^m}\om_h^m= \frac {m!}{\pi^m} \,d\vol_M\,.$$
\end{theo}
 
\subsection{Berman's canonical \kahler point process}

In a series of articles, \cite{Ber11,Ber14,Ber17,Ber18,Ber20,Ber21}, R. Berman investigated determinantal
 point processes on  K\"ahler manifolds defined in terms of Bergman kernels
and related geometric invariants. 

Let $M$ be a compact \kahler manifold and suppose that the canonical line bundle $K_M=\bigwedge_mT^*_M$ is ample. It was shown by Aubin \cite{Aubin} and by Yau \cite{Yau} that $M$ carries a K\"ahler--Einstein metric $\om_{KE}$; i.e., the metric has constant scalar curvature: Ric$\,\om_{KE}=-\om_{KE}$.  
Berman constructs  determinantal point processes $[{\bf z}_k]$ and uses their  expected empirical measures to obtain \kahler metrics on $M$ converging to $\om_{KE}$  as $k\to\infty$. 

Whereas the zero point process above depends on the choice of the metric (i.e., on the measure $\gamma_{h^k,dV}$), Berman's canonical point processes are   independent of the choice of metric. As before (with $L=K_M$), let
$\{S^k_1, \dots,S^k_{d_k}\}$ be a basis for the pluricanonical system $H^0(M,K_M^k)$, where $d_k=\dim H^0(M,K_M^k)$.
 Berman's canonical probability measure (or point process) 
on the configuration space $M^{d_k}$ of $d_k$ points in $M$ is the probability measure $\mu_{k}^B$ defined by
$$\mu^B_k : = \frac{1}{\la_k} \left| \det S^{(k)} (z_1, \dots, z_{d_k}) \right|^{\frac{2}{k}},$$
where $\la_k$ is a normalizing constant and $$\det S^{(k)} (z_1, \dots, z_{d_k}): = \det \begin{pmatrix} S_j^k(z_l) \end{pmatrix}_{1\le j,l\le d_k}$$ is a holomorphic section of the pluricanonical bundle
$ K^k_{M^{ d_k}} \to M^{d_k}=M\times \cdots\times M$. Thus $|\det S^{(k)}|^{2/k}$ is a semi-positive section of $K_{M^{ d_k}}\wedge \overline{K_{M^{ d_k}}}$, i.e. a positive measure on $M^{ d_k}$.
Changing the basis changes $|\det S^{(k)}|^{2/k}$ by a constant factor, and we divide by $\la_k$ so that $\mu^B_k$ is a well-defined probability measure.

As mentioned at the beginning of Section \ref{section-zeros}, a point ${\bf z}_k=(z_1,\dots,z_{d_k})\in M^{d_k}$ gives rise to the empirical probability measure on $M$, $$[{\bf z}_k]:= \frac 1{d_k}\sum_{j=1}^{d_k}\delta_{z_j}\,.$$ Thus (after dividing out by the symmetric group) we can consider $\mu^B_k$ to be a probability measure on the space of discrete  probability measures on $M$.  

Berman then obtains canonical sequences of \kahler forms and volume forms on $M$ converging to the K\"ahler-Einstein metric and  volume, respectively:
\begin{theo}\label{Bvolume-limit}\cite{Ber17} Let $M$ be a compact \kahler manifold such that $K_M$ is ample.  \begin{enumerate}

\item[i)] \hspace{1.5in}$\E_{\mu^B_k}[{\bf z}_k] \to c\,\om_{KE}^m\,, \quad \mbox{as }\ k\to\infty,$\\[4pt] where $c$ is a normalizing constant;\\
 \item[ii)] Writing $\E_{\mu^B_k}[{\bf z}_k]= f_k \,(i^{m^2}\eta\wedge\bar\eta)$ over an open set $U$, where $\eta$ is a nonvanishing holomorphic $m$-form, the \kahler form 
 $$\om_k:= \frac i{2\pi}\ddbar \log f_k\to \om_{KE} \,, \quad \mbox{as }\ k\to\infty.$$
 \end{enumerate}
\end{theo}
Note that the \kahler forms $\om_k$ are globally defined and independent of the choice of $\eta$.  Part (i) is an analog of Theorem \ref{volume-limit}. The constant $c$  is chosen so that $\int_Mc\,\om_{KE}^m= c\,{m!}\vol_{KE}(M)=1$.

\section{\label{RanBERG} Random Bergman metrics}

In this section, we discuss a recent direction to stochastic K\"ahler geometry: the study  of random K\"ahler metrics in a fixed class
$\kcalomega$  and their approximations by random Bergman metrics, as given by Ferrari, Klevtsov, and Zelditch \cite{FKZ12a,FKZ12b,KZ16}.  Here, $\om_0=\om_h=\pi c_1(L,h)$ is the \kahler metric of a positive line bundle $(L,h)\to (M,\om_h)$ and $\kcalomega$ is the infinite dimensional space of 
 K\"ahler metrics    $\omega\in[\omega_0]$, the cohomology class of $\omega_0$. The space of all K\"ahler metrics $\kcalomega$
on $M$ in the K\"ahler class $[\omega_0]$ is parametrized as 
\begin{equation}\label{KCAL} 
\mathcal K_{[\omega_0]}= \{\phi\in C^{\infty}(M)/\mathbb R\,:\, \omega_0+i\p\bp\phi>0\}.
\end{equation} 
The motivation to study rather general types of random K\"ahler metrics originates in some sense in Polyakov's approach to  quantum gravity. In complex dimension one, it has led to an explosion of articles on LQG (Liouville quantum gravity). In keeping with our emphasis on higher dimensional K\"ahler manifolds, we do not review the voluminous literature on LQG but only the random K\"ahler metrics studied in \cite{FKZ12a,FKZ12b,KZ16}.
To endow \eqref{KCAL} with an interesting probability measure  is very difficult because of its infinite dimensionality. In LQG, one specific measure is studied and it is induced by a well-studied Gaussian field, the Gaussian free field. More precisely, it is a renormalized version of the exponential of the GFF and is known as the Gaussian multiplicative chaos. In higher dimensions, there is no parallel construction and one has to start from scratch. The main idea is to define a sequence $\mu_k$ of  probability measures on finite dimensional spaces $\bcal_k$ of Bergman metrics and then to study their 
 limits.

We begin by describing the spaces  $\bcal_k$ and then focus on one specific choice of probability measure induced by Brownian motion on $\bcal_k$ with respect to its symmetric space Riemannian metric. The space $\bcal_k$ of Bergman metrics of degree $k$ is the space of metrics given by the pullbacks of Fubini-Study metrics by the Kodaira map for $H^0(M,L^k)$. I.e., let $\{\sigma_1,\dots,\sigma_{d_k}\}$ be a basis for $H^0(M,L^k)$, and let \beq\label{iota}\iota_{\sigma}  = [\sigma_1, \dots, \sigma_{d_k}] :M \to \CP^{d_k-1}\,.\eeq
Since positive line bundles are ample, we can choose $k$ sufficiently large so that \eqref{iota} is an imbedding (for all bases $\{\sigma_j\}$ of $H^0(M,L^k)$).  The associated Bergman metric is
\beq\label{Bergmet}  \frac 1k\iota_\sigma^* \omega_{FS} = \frac{i}{2k} \, \ddbar \log \sum_{j = 1}^{d_k} |f_j|^2\,,\eeq  where $\sigma_j=f_je_L^{\otimes k}$ for a local frame $e_L$. The space $\bcal_k$ of Bergman metrics of degree $k$  then consists of all metrics of the form \eqref{Bergmet}. 

The space $\bcal_k$ can be parametrized by the symmetric space $\SU(d_k)\backslash\SL(d_k,\C)$ as follows: let $\{S_1^k,\dots,S_{d_k}^k\}$ be a fixed  orthonormal basis with respect to the inner product \eqref{inner} induced by  $h$ and $\om_0$, and write $S^k_j=F^k_je_L^{\otimes k}$ as above. For matrices $A=(A_{jl})\in\SL(d_k,\C)$, we let $\sigma^A_j=\sum_lA_{jl}S^k_l$.    Then $\sigma^A=\{\sigma^A_1,\dots \sigma^A_{d_k}\}$ is a basis for $H^0(M,L^k)$, and the associated Bergman metric is
\beq\label{BergmetP}\frac 1k(\iota_{\sigma^A})^*\om_{FS} =\frac{i}{2k} \,  \ddbar \log \sum_{jl}\bar F^k_jP_{jl}F^k_l= \om_0+\frac{i}{2k} \,  \ddbar \log \sum_{jl}\left\|\bar S^k_jP_{jl}S^k_l\right\|_{h^k}\,,\eeq
where   $P=A^*A$ is in the space of positive definite Hermitian $d_k\times d_k$ matrices with determinant one. We denote this space by $\pcal_{d_k}$ and note that $\SU(d_k)\backslash\SL(d_k,\C) \cong \pcal_{d_k}$ via the map $A\mapsto A^*A$.

For matrices $P=A^*A\in\pcal_{d_k}$, we define the {\it Bergman potential}  \beq\label{phiP}\phi_P:=\frac 1{2k} \,  \log \sum_{jl}\bar F^k_jP_{jl}F^k_l\,,\eeq and we let \beq\label{omegaP} \om_P:=\frac 1{2k}(\iota_{\sigma^A})^*\om_{FS} =i\,\ddbar\phi_P\eeq denote  the corresponding Bergman metric.\footnote{Here, we are using the convention that $\pi \om_{\FS}$ is in the Chern class of the hyperplane section bundle $\ocal(1)\to\CP^m$, and thus $[\pi\om_P]=c_1(M,L)$ for $\om_P\in\bcal_k$. Hence in this article, $\om_P$ and $\phi_P$ equal $\frac12$  the corresponding terms in \cite{FKZ12b, KZ16}.\label{f1}} In particular, $\phi_{\I_k}=\phi_h+\frac 1k\log \|B_k(z,z)\|$, and
\beq\label{omegaI}\om_{\I_k}=\om_0+\frac i{2k}\ddbar \log\| B_k(z,z)\|= \om_0+O\left(\frac 1{k^2}\right),\eeq by \eqref{seg}, where  $\I_k$ is the $d_k\times d_k$ identity matrix.

\subsection{Heat kernel measures}

Given an orthonormal basis $\{S_j^k\}$ of $H^0(M,L^k)$ with respect to $h$ and $\om_0=\om_h$, the space $\bcal_k$ can be identified with the  symmetric space $\pcal_{d_k}$ via equations \eqref{BergmetP}--\eqref{omegaP}.  The
 general   question is to find
 sequences $\{d\mu_k\}$ of measures on $\bcal_k$ that are independent of the choices of
the basis $\{S^k_j\}$ and which vary in a simple way under the change of the reference point $\omega_0\in\kcalomega$ and have good asymptotic properties
as $k \to \infty$. Such measures can be given as
the  heat kernel  measures
 \begin{equation} \label{HEAT} d\mu_k^t(P): = p_k(t, \I_k, P)\, dP, \end{equation} 
 where $dP$ is Haar measure on $\pcal_{d_k}$, and  $p_k(t, P_1, P_2)$ is the heat kernel of the symmetric space $\pcal_{d_k}$.  The measure is invariant under the action of the unitary group ${\rm U}(d_k)$ and thus is independent of the choice of the orthonormal basis of sections $\{S^k_j\}$ used for the  matrix-metric identification in \eqref{BergmetP}.   Then \eqref{HEAT} is the probability measure on $\bcal_k$ induced by 
 Brownian motion on $\pcal_{d_k}$ starting at the identity $\I_k$ at time $t=0$.

In this section we review results of \cite{KZ16} on  the behavior of the heat
kernel measure \eqref{HEAT} on $\bcal_k$ as $k \to \infty$.  The heat kernel measure is only one among many possible measures to study; we choose it because it has a simple geometric and probabilistic  interpretation and because we obtain surprisingly explicit formulae for its 
correlation function.
However, it is so closely tied to the symmetric space geometry of positive Hermitian matrices that it does not reflect the deeper geometric aspects of
$\bcal_k$. At the end of this section, we propose a model which does go deeper, namely the Calabi metric measure on $\bcal_k$. However it is difficult to obtain analytic expressions for the key probabilistic objects for this Calabi model.

It was shown in \cite{FKZ12b}, that for all probability measures $\nu$ on $\pcal_{d_k}$ on $\pcal_{d_k}$ invariant under the ${\rm U}(d_k)$ action $P\mapsto U^*PU$, one has 
 \beq\label{1ptphi} \E_\nu\phi_P=\phi_{\I_k}=\frac 1{2k}\log\|B_k(z,z)\| - \log\|e_L(z)\|\,,\eeq and thus by \eqref{seg}
 \beq\label{1pt} \E_\nu\om_P=\om_{\I_k}=\om_0+O\left(\frac 1{k^2}\right).\eeq In particular \eqref{1pt} holds for the heat kernel measures.

However, the two-point correlations depend on the choice of invariant measure. The two-point correlations for the heat kernel measures \eqref{HEAT} were given in \cite{KZ16}, where it was shown that the correlations have the form  
\beq\label{final2p} 
\E_{\mu^t_k}\, \phi_P(z)\phi_P(w)=\phi_{\I_k}(z)\phi_{\I_k}(w)+\frac1{4k^2} I_{2}(t,\beta_k(z,w)),
\end{equation}
 where \begin{equation}\label{rhointro} \beta_k(z,w) = \frac{\|B_k(z, w)\|^2}{\|B_k(z, z)\|\,\| B_k(w, w)\|} = P_k(z,w)^2 \end{equation}
is the {\it Berezin kernel}, and \beq\label{I2k}
\frac\d{\d x}\, I_{2}(t, x)
=\frac{2t}x-\frac{e^{-t/2}}{\sqrt{2\pi t}}\frac{\sqrt{1-x}}{x}\int_{-\infty}^\infty d\lambda\, \frac{\,e^{-\frac1{2t}\lambda^2}\cosh\lambda}{\sqrt{\coth^2\lambda-x}}\log\frac{\sqrt{\coth^2\lambda-x}+\sqrt{1-x}}{\sqrt{\coth^2\lambda-x}-\sqrt{1-x}}.\eeq 
 
It follows from \eqref{1ptphi}  and \eqref{final2p} that
\beq \label{varphi}\Var(\phi_P) = \frac1{4k^2} I_{2}(t,\beta_k)\,,\eeq
where $\Var=\Var_{\mu^t_k}$ is given by Definition \ref{varcurrent}. Note that $$\big(\Var(\phi_P)\big)(z,w)={\rm Cov}\big(\phi_P(z),\phi_P(w)\big)\,.$$
Furthermore, differentiating \eqref{varphi}, we have
\begin{equation}\label{vc2}
{\bf Var}\big(\omega_{P}\big)=\frac1{4k^2}
 (i\ddbar)_z\,(i\ddbar)_w \,I_{2} (t,\beta_k(z,w))\,.\eeq
Formula \eqref{vc2} says that $I_{2}(t,\beta_k)$  is the pluri-bipotential of the variance of the \kahler metric for the heat kernel measure at time $t$. In the Riemann surface case ($\dim M=1$), the variance of the area of a domain $U\subset M$ is given by
$$\var\left(\int_U\om_P\right)=\int_{U\times U}\Var(\om_P)= -\frac{1}{4k^2}\int_{\d U\times\d U} \d_z\,\d_w \,I_{2} (t,\beta_k(z,w))\,.$$

If we fix $k$ and let $t\to\infty$. Then \beq\label{ktoinfinity}\frac\d{\d x}\, I_{2}(\infty, x):=\lim_{k \to \infty} \frac\d{\d x}\, I_{2}(t, x)  = -\frac{\log (1 - x)}{x}\,,\eeq (see \cite{KZ16}) and thus $$I_2(\infty,x)=\Li_2(x)\,.$$
Therefore for fixed $k$, the variance of random \kahler metrics with respect to the heat kernel measure $\mu^t_k$  on the space $\bcal_k\cong \pcal_{d_k}$ converges to the variance of the scaled zero current $\frac{\pi}{k}Z_{s^k}=\frac ik\ddbar\log|f|$ of a random section $s^k=fe_L^{\otimes k}\in H^0(M,L^k)$ (given in Theorem~\ref{BIPOT}) as $t\to\infty$.  In fact, in the $t\to\infty$ limit, the  random metrics converge to random zero divisors regarded as singular metrics, in the sense described in \cite[\S5.1]{KZ16}.

Since \eqref{vc2} gives an exact formula for any $t, k$, one may also consider a variety of limits as $t \to \infty, k \to \infty$ in some relation.
There is a natural choice of coupled limit motivated by the metric asymptotics of $\bcal_k$. 
If we  rescale
the Cartan-Killing (CK) metric  $g_{CK,k}$ on $\pcal_{d_k}\approx \bcal_k\subset \kcalomega$ as $ g_k =  \epsilon_k^2 g_{CK, k}$, with
  $\epsilon_k = k^{-1} d_k^{-1/2}$, then $g_k \to g_M$ on $T \bcal_k$. Here $g_M$ is
 the Mabuchi metric  on $\kcalomega$, i.e. the  Riemannian metric on $\kcalomega$ 
defined by $||\delta \phi||^2_{\phi_0} = \int_M (\delta \phi)^2 \omega_{\phi}^m/m!
$, where $\omega_{\phi} = \omega_0 + i \ddbar \phi$.  Thus, a ball of radius one with respect to the usual CK metric $g_{CK,k}$ has radius
approximately $\epsilon_k$ with respect to the Mabuchi distance. 
With the rescaling $g_k =\epsilon_k^2 g_{CK, k}$, the corresponding
Laplacian scales as $\Delta_{g_k} \mapsto \epsilon_k^{-2} \Delta_{g_{CK, k}}$. It follows
that the heat operator scales as
$$\exp t \Delta_{g_k} = \exp t \epsilon_k^{-2} \Delta_{g_{CK, k}}. $$
In effect, it is only the time that is rescaled, and the rescaled heat kernel
is $p_k(\epsilon_k^{-2} t, \I_k, P). $

We therefore study the  metric scaling limit with $t \mapsto t \epsilon_k^{-2}$ and
 evaluate $I_{2}( \epsilon_k^{-2} t, \beta_k)$ asymptotically as $k \to \infty$. This scaling 
keeps the $d_k$-balls of uniform size as $k \to \infty$ with respect to the limit
Mabuchi metric. Thus, as $k$ changes, the Brownian motion with respect to
$g_k$ probes distances of size $t$ from the initial metric $\omega_0$ for all $k$.
It turns out that  \beq\lim_{k \to \infty}  I_{2}\left(\epsilon_k^{-2} t, \beta_k(z,z+u/\sqrt k) \right) =\Li_2(e^{-|u|^2})\,.\eeq

\subsection{The Calabi model}

The Calabi metric is the natural (background independent)
$\lcal^2$ metric on either $\kcalomega$ or $\bcal_k$. If $ \dot{\omega} = \delta \omega \in T_{\omega} \kcalomega$,
$$
|| \delta \omega ||_{C}^2= \int_M ||\delta \omega(z)||^2_{\omega}
dV_{\omega}. $$
It is the restriction to a \kahler class  of the deWitt-Ebin metric on metric
tensors \cite{E70,DeWitt}.
In terms of relative \kahler  potentials, $ \dot{\omega} =  i\ddbar \dot{\phi}$,  the Calabi metric inner product  is,
\begin{equation} \label{CM} 
||\Delta_{\omega} \dot{\phi} ||^2_{C, \omega} = \int_M | \Delta_{\omega} \dot{\phi} |^2
\,dV_\omega.  \end{equation}
It is known that the sectional curvatures of $(\kcalomega, g_C)$ are all
equal to $1$, i.e. this Riemannian manifold is an open subset of the infinite dimensional sphere of constant curvature $1$ (see \cite{Calam12}). 
The  finite dimensional approximations to the Calabi metric should approximate domains in  finite
dimensional spheres.  Hence, 
 
 \begin{conj} \label{FINITE} The Calabi volume $\rm{Vol}_k (\bcal_k)$ with 
 respect to $G |_{\bcal_k}$  is finite for each $k$. 
 \end{conj}

If the conjecture is true, one obtains a purely geometric sequence of probability measures. This could give a rigorous definition of the Polyakov path integral over metrics, which used the Calabi metric to define its volume form. Polyakov also used a power of the determinant of the Laplacian, which could also be implemented in the Bergman approximation.

\end{document}